\documentclass{amsart}
\usepackage[toc,page]{appendix}
\usepackage{amsmath}
\usepackage{amssymb}
\usepackage{algorithm}
\usepackage{algpseudocode}
\usepackage{amsthm}
\usepackage{graphicx}
\usepackage{subcaption}
\usepackage{verbatim}
\usepackage{epsfig}
\usepackage{hyperref}
\usepackage{float}
\usepackage{color}
\usepackage{fullpage}
\usepackage{amsfonts}
\usepackage{amscd}
\usepackage{amssymb}
\usepackage{bbm}
\usepackage{graphics}
\usepackage{amsmath, bm}
\usepackage{amsbsy}
\usepackage[english]{babel}
\usepackage{bbm}
\usepackage{bm}
\usepackage{dsfont}
\usepackage{mathrsfs}
\usepackage{color}
\usepackage{enumitem}
\usepackage[table,xcdraw]{xcolor}  
\usepackage{cleveref}
\usepackage[toc,page]{appendix}
\DeclareFontFamily{OML}{rsfs}{\skewchar\font'177}
\DeclareFontShape{OML}{rsfs}{m}{n}{ <5> <6> rsfs5 <7> <8> <9>
rsfs7 <10> <10.95> <12> <14.4> <17.28> <20.74> <24.88> rsfs10 }{}
\DeclareMathAlphabet{\mathfs}{OML}{rsfs}{m}{n}
\newcommand{\hldots}{\raisebox{0.5ex}{$\ldots$}}
\newcommand{\ind}{\mathbf{1}}

\newcommand{\bae}{\begin{equation}\begin{aligned}}
\newcommand{\eae}{\end{aligned}\end{equation}}

\newtheorem{theorem}{Theorem}[section]

\theoremstyle{definition}
\newtheorem{example}{Example}

\usepackage[margin=1in]{geometry}
\usepackage{amsmath,amssymb,amsthm,mathtools}
\usepackage{bm}
\usepackage{booktabs}
\usepackage{graphicx}
\usepackage{enumitem}
\usepackage{microtype}
\usepackage{algorithm}
\usepackage{algpseudocode}
\usepackage{tikz}
\usetikzlibrary{arrows.meta,positioning,shapes.geometric,fit,calc,decorations.pathmorphing}

\makeatletter
\newcommand{\liftheadings}{%
  \renewcommand{\subsection}{\@startsection{subsection}{2}%
    \z@{.7\linespacing\@plus\linespacing}{.5\linespacing}%
    {\normalfont\scshape\centering}}%
  \renewcommand{\subsubsection}{\@startsection{subsubsection}{3}%
    \z@{.5\linespacing\@plus.7\linespacing}{-.5em}%
    {\normalfont\bfseries}}%
  \renewcommand{\paragraph}{\@startsection{paragraph}{4}%
    \z@{.5\linespacing\@plus.7\linespacing}{-.5em}%
    {\normalfont\itshape}}%
}
\makeatother

\begin{document}
\numberwithin{equation}{section} \numberwithin{figure}{section}
\title{Global genetic inference at the scale of mass disasters}
\author{Youri Weesie, Ritsaart Bergsma, Aafko Boonstra}

\begin{abstract}
DNA-based disaster victim identification is commonly performed sequentially, with matches accepted before all genetic evidence is considered jointly. However, early matches may be incorrect or inconsistent with the globally most likely assignment, and treating them as final can affect subsequent identifications. Joint analysis avoids this dependence on decision order, but combinatorial growth has historically made it prohibitive in large cases. In this paper, we show that the general joint identification problem is NP-hard and develop a local-to-global framework for joint genetic inference at scale. We reduce the search using genetic compatibility tests, likelihood-based pruning of unlikely assignments, and symmetry reductions. We represent pedigree likelihoods as polynomials and reformulate the Elston–Stewart algorithm for efficient parallel likelihood evaluation on GPUs. Together, these advances yield a manageable set of candidate assignments and their likelihoods, from which we construct an integer linear program that resolves competing identifications and determines the global maximum-likelihood assignment. In small benchmarks, the method recovers the exhaustive-search optimum with up to \(99.96\%\) fewer likelihood evaluations. We also apply it to simulated disasters involving up to 10 thousand victims and missing persons. The framework is not limited to a single global optimum but yields a list of probable alternative assignments, thus providing a flexible tool to be used alongside other forensic evidence.
\end{abstract}
\maketitle
\section{Introduction}
DNA-based \textit{disaster victim identification} (DVI) is a well-established field in forensic science, with notable examples including the 2026 Nepal floods \cite{nepal2026}, the 2026 New Year's Eve bar fire in Crans-Montana, Switzerland \cite{crans2026}, the Malaysia Airlines Flight MH17 \cite{vermeij2022} and the 2004 Thailand tsunami \cite{brenner2006}.
While various sources of evidence can support victim identification, this paper only considers genetic evidence. In a DVI context, this evidence is usually split into two categories: the post-mortem (PM) DNA obtained from victim material and ante-mortem (AM) DNA provided by biological relatives of missing persons, optionally supplemented with auxiliary information such as sex.  This paper is concerned with the computational aspects of identification cases involving multiple and possibly related victims, after the data have been collected \cite{prinz2007, sweet2010}.

A standard approach matches a victim to a missing person when that specific match is strongly supported by the genetic evidence. When the identified individual belongs to a reference family, the victim’s profile can then be incorporated into the reference data to help identify other missing relatives in the same family. This process is repeated until no new identifications can be made. Following \cite{vigeland2021, DVIstat}, we refer to this as a \textit{sequential method}.

Sequential methods, such as those implemented in Bonaparte \cite{vandongen2011, slooten2010} and Familias \cite{egeland2000, kling2017} can perform well when individual identifications are strongly supported and hence effectively independent. However, each decision is made using only part of the evidence available for the complete DVI case. Once an assignment has been accepted, it is generally treated as fixed, and later decisions are conditioned on that assignment rather than being obtained by reconsidering all assignments simultaneously. An early incorrect assignment can therefore propagate through the remaining analysis and produce additional incorrect or order-dependent identifications \cite{vigeland2021}. Sequential methods may also miss assignments that receive strong support jointly but only weak support when evaluated pairwise \cite{DVIstat}.

\textit{Joint identification}, as advocated in \cite{vigeland2021, DVIstat}, addresses these limitations. This approach seeks assignments that maximize the joint likelihood over all identifications. To the best of our knowledge, the $R$-package \texttt{dvir} \cite{vigeland2021book} is the only freely available software that supports joint DVI computations. Although it contains several elements to speed up the evaluation process, it is limited in the size of the problems it can consider. In this work, we contribute to the literature on joint DVI by combining local inference within reference families with global optimization. We reduce the assignment space through genetic compatibility filtering, likelihood-based pruning, and symmetry reduction, and use a polynomial reformulation of the Elston–Stewart algorithm for efficient parallel likelihood evaluation on GPUs. An integer linear program then combines the resulting assignments and likelihoods to determine the global maximum-likelihood assignment. The source code of our new approach is freely available from \url{https://github.com/[]}.

In \Cref{methods}, we describe the core components of our approach for DNA-based DVI problems, after which we present the results and the method's limitations in \Cref{results} and \Cref{discussion} respectively. Before presenting the method, we emphasize that DNA evidence is interpreted probabilistically. With sufficient genetic data, conclusions can often be drawn with high confidence; however, the most likely explanation may nevertheless differ from the ground truth. This is a natural limitation of DNA; the probabilistic uncertainty cannot be eliminated entirely, but this uncertainty can be reduced by allowing additional, potentially non-genetic information to be incorporated \cite{vergeer2023, boonstra2025}. For this reason, the method returns a ranked list of the most likely configurations rather than returning only the single best explanation, so that they can be evaluated alongside the other case-specific evidence. For forensic practitioners, such a list generally provides an adequate basis for subsequent forensic decision-making \cite{DVIstat}. 

\section{Methods}\label{methods}
\subsection{Definitions} We begin with a formal outline of the DVI problem. Following \cite{vigeland2021}, the first component is the \emph{ante-mortem} (AM) data. This is a collection
\begin{equation}
    \mathcal{M} = \{M_1,\ldots,M_t, O_1, \dots, O_{s'}\},
\end{equation}
of $t$ missing persons accompanied by $s'\leq s$ \textit{open spots}. This number $s'$ is bounded by the \emph{post-mortem} (PM) data, written as
\begin{equation}
    \mathcal{V} = \{V_1,\ldots,V_s, *_1, \dots, *_{t'}\},
\end{equation}
which consists of the genetic profiles of $s$ victims together with $t' \leq t$ \textit{untyped} `dummy' victims. The open spots and untyped victims are created to formally allow the option for victims and missing persons to remain unidentified. We require that $s+t'=t+s'$ and call a problem \textit{fully open} if both $s=s'$ and $t=t'$. 

In addition to the AM data, we have (disjoint) reference families $F_1,\ldots,F_J$. We denote $I_j \subseteq \mathcal{M}$ for the missing persons of family $F_j$. To accommodate for the open spots, we define $I_0 := \bigcup_{i=1}^{s'}\{O_i\},$ so that $\mathcal{M}=\bigcup_{j=0}^J I_j.$

The individuals in a reference family are of the following three types:
\begin{itemize}
    \item \textit{Known members}: relatives for whom DNA profiles are available. These supply AM genetic information used to evaluate the plausibility of each PM--missing match.
    
    \item \textit{Untyped members}: individuals required for the pedigree structure but without observed genotypes. Their genetic states are unknown and must be summed or marginalized over when computing likelihoods.
    
    \item \textit{Missing members}: individuals whose identities are sought. They correspond to the elements of $\mathcal{M}$ and must be matched to PM samples if genetically supported.
\end{itemize}

For a graphical representation of these concepts, see e.g. \Cref{runex} below. Next, we represent a joint identification by an \emph{assignment}. This is a map
\begin{equation}
   x : \mathcal{M} \to \mathcal{V} , 
\end{equation}
where $x(M_j) = V_i$ indicates that PM sample $V_i$ is identified as missing person $M_j$. In this notation, $x(M_j)=\ast_i$ indicates that $M_j$ has not (yet) been identified, and similarly, $x(O_i)=V_j$ in turn reflects that victim $V_j$ is not matched to any reference family. Assignments must of course be bijective as a victim should not be identified with more than one missing person.

Even for moderate values of $s$ and $t$, the number of possible assignments can be extremely large \cite{vigeland2021}. Exhaustive evaluation of all possible $x$ to find the most likely assignment is therefore computationally prohibitive. The following theorem makes this precise. 

\begin{theorem}
    The DVI problem is NP-hard.
\end{theorem}

We refer to \Cref{npbewijs} for a formal statement and proof.

A key observation towards an optimal assignment is that the dependencies that drive the complexity occur locally, i.e. only within reference families. The likelihood of a global assignment therefore decomposes into a product of conditional likelihoods per family. As we will see, calculating all configurations per family is still not an easy task, but already much better than counting all possible global assignments.

Once this is done, we rephrase the DVI problem into the language of \textit{integer linear programs} (ILPs) \cite{williams2013}. The resulting ILP is very similar to an \emph{exact weighted set cover problem} \cite{balas1976}, which can usually be solved with the help of standard ILP solvers \cite{song2026}. This process returns the globally optimal assignment and, if desired, a ranked list of near-optimal alternatives.

In the remainder of this section, we will cover the core principles of the implementation. 

\subsection{ILP formulation} In order to rewrite a DVI problem into an ILP, we introduce some notation. We start by writing $d_j$ for the number of missing persons in family $F_j$. We let $\mathcal{V}^{d_j}$ be the tuples of cardinality $d_j$ of $\mathcal{V}$ representing the space of all possible assignments for the missing persons, in a family $F_j$. We furthermore assume an ordering on $\mathcal{V}^{d_j}$ and write $s_{i,j} \in \mathcal{V}^{d_j}$ for its $i$-th element. 

Then, for any assignment $s_{i,j} \in \mathcal{V}^{d_j},$ we denote the conditional likelihood with
\begin{equation}\label{ll}
    P(s_{i,j} \mid F_j)
\end{equation}
for each $s_{i,j}$. In this work, the likelihood \cite{elston1971, lange2002} is calculated assuming Mendelian segregation, unlinked markers, linkage equilibrium and mutation models.
We further define the \textit{weights}
\begin{equation}
  w_{i,j} := \log (P(s_{i,j} \mid F_j))
\end{equation}
as the corresponding natural log-likelihood of the victims placed in $F_j$, and define $x_{i,j}$ as an indicator function that is 1 if the restriction $x_{|I_j}$ is equal to $s_{i,j}$. 

For the global log-likelihood $L(x)$ of the assignment $x,$ we can now write

\begin{equation}
  L(x) := \sum_{j=0}^J \sum_{i=1}^{n_j} w_{i,j}x_{i,j},
\end{equation}
where $n_j$ denotes the number of candidate assignments to $F_j$.
The objective is to identify the assignment $x$ that maximizes the likelihood $L(x)$ among all feasible mappings, given the PM and AM data. With the above notation, we can reformulate this as

\begin{equation}\label{obj}
\max \sum_{j=0}^J \sum_{i=1}^{n_j} w_{i,j}x_{i,j}
\end{equation}
subject to
\begin{equation}\label{eis1}
\sum_{i=1}^{n_j} x_{i,j} = 1 \textup{ for all } 0 \leq j \leq J,
\end{equation}
meaning that every family gets one single assignment for the missing persons and
\begin{equation}\label{eis2}
\sum_{j=0}^J \sum_{i: v \in s_{i,j}} x_{i,j} = 1 \textup{ for all } v \in \mathcal{V},
\end{equation}
meaning that every victim is only assigned once. In \eqref{eis2}, $v \in s_{i,j}$ indicates with slight abuse of notation that $v$ occurs as an entry in $s_{i,j}$. Note that \eqref{eis1} and \eqref{eis2} effectively require $x$ to be an injection and thus a bijection from $\mathcal{M}$ to $\mathcal{V}$.

The above reduces the problem to finding the values $w_{i,j}$ for each $j$. That is, for each family $F_j$, we need to find the likelihood of all possible assignments. Generally, this still poses computational problems because the number of all possible assignments quickly explodes as $s$ and $t$ become larger, see \cite{vigeland2021} for an exact count.

Before we discuss the main strategies to overcome this hurdle, we present the following toy DVI example of manageable size. 

\begin{example}
Consider \Cref{runex}. In this example, the PM data consists of four victims $V_1$ through $V_4$ that have been recovered, each with a known DNA profile. The AM data involves two families $F_1$ and $F_2$. Both contain typed reference relatives and a number of missing individuals whose identities are sought. Family $F_1$ has two missing persons $M_1$ and $M_2$, with $M_1$ being a child of the typed reference $R_1$. Family $F_2$ has two missing persons $M_3$ and $M_4$, whose parents are the typed references $R_2$ and $R_3$. The problem is fully open, as we have 4 open spots ($F_0$) and 4 dummy victims $*_1$ through $*_4$.

\tikzset{
  ped male/.style   = {draw, rectangle, minimum size=5.5mm, inner sep=0pt},
  ped female/.style = {draw, circle,    minimum size=5.5mm, inner sep=0pt},
  ped ref/.style    = {fill=black!72, text=white, thick},
  ped miss/.style   = {dashed, thick},
  ped unk/.style    = {fill=black!7},
  ped vic/.style    = {fill=black!0},
}
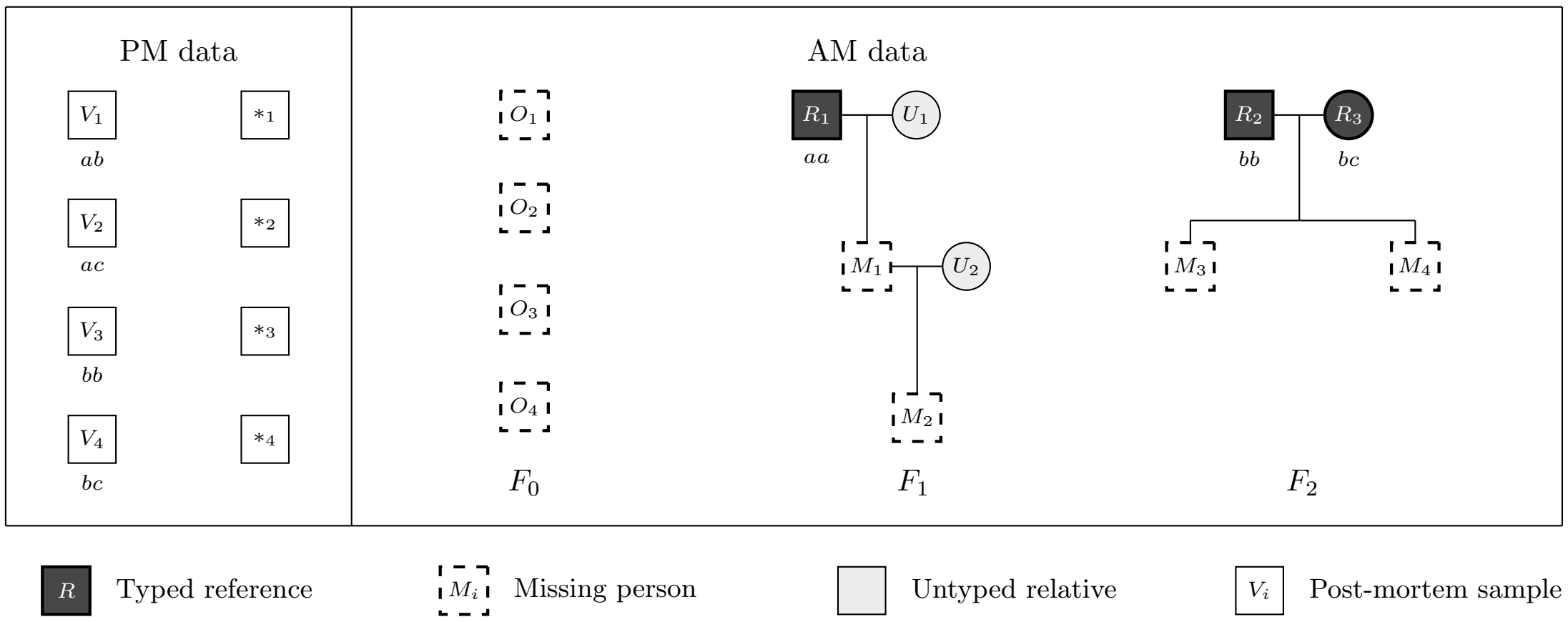
\begin{figure}[H]
\centering
\begin{tikzpicture}[scale=0.8, every node/.style={transform shape}]

\draw (-9,2.00) -- (9,2.00);

\draw (-9,-4.00) -- (9,-4.00);

\draw (9,2.00) -- (9,-4.00);

\draw (-9,2.00) -- (-9,-4.00);

\node at (-7,1.5) {PM data};

\node[ped male, ped vic] at (-8,0.75) {\scriptsize $V_1$};
\node[ped male, ped vic] at (-8,-0.50) {\scriptsize $V_2$};
\node[ped male, ped vic] at (-8,-1.75) {\scriptsize $V_3$};
\node[ped male, ped vic] at (-8,-3) {\scriptsize $V_4$};

\node[ped male, ped vic] at (-6,0.75) {\scriptsize $*_1$};
\node[ped male, ped vic] at (-6,-0.50) {\scriptsize $*_2$};
\node[ped male, ped vic] at (-6,-1.75) {\scriptsize $*_3$};
\node[ped male, ped vic] at (-6,-3) {\scriptsize $*_4$};

\node at (-8,0.25) {\scriptsize $ab$};
\node at (-8, -1) {\scriptsize $ac$};
\node at (-8,-2.25) {\scriptsize $bb$};
\node at (-8,-3.5) {\scriptsize $bc$};

\draw (-5,2.00) -- (-5,-4.00);


\node[ped male, ped miss] at (-3,0.75) {\scriptsize $O_1$};
\node[ped male, ped miss] at (-3,-0.325) {\scriptsize $O_2$};
\node[ped male, ped miss] at (-3,-1.5) {\scriptsize $O_3$};
\node[ped male, ped miss] at (-3,-2.625) {\scriptsize $O_4$};


\node at (0.96,1.5) {AM data};
\draw (0.38,0.75) -- (1.53,0.75);
\draw (0.96,0.75) -- (0.96,-0.71);
\draw (1.22,-1) -- (2.17,-1);
\draw (1.54,-1) -- (1.54,-2.49);
\node[ped male, ped ref] at (0.38,0.75) {\scriptsize $R_1$}; \node at (0.38, 0.25) {\scriptsize $aa$};
\node[ped female, ped unk] at (1.53,0.75) {\scriptsize $U_1$};
\node[ped male, ped miss] at (0.96,-1) {\scriptsize $M_1$};
\node[ped female, ped unk] at (2.11,-1) {\scriptsize $U_2$};
\node[ped male, ped miss] at (1.54,-2.75) {\scriptsize $M_2$};

\node at (-3,-3.5) {$F_0$};
\node at (1.5,-3.5) {$F_1$};
\node at (6,-3.5) {$F_2$};

\draw (5.38,0.75) -- (6.53,0.75);
\draw (5.96,0.75) -- (5.96,-0.47);
\draw (4.70,-0.47) -- (7.30,-0.47);
\draw (4.70,-0.47) -- (4.70,-0.71);
\draw (7.30,-0.47) -- (7.30,-0.71);
\node[ped male, ped ref] at (5.38,0.75) {\scriptsize $R_2$}; \node at (5.38,0.25) {\scriptsize $bb$};
\node[ped female, ped ref] at (6.53,0.75) {\scriptsize $R_3$}; \node at (6.53,0.25) {\scriptsize $bc$};
\node[ped male, ped miss] at (4.70,-1) {\scriptsize $M_3$};
\node[ped male, ped miss] at (7.30,-1) {\scriptsize $M_4$};

\node[ped male, ped ref] at (-8.3,-4.75) {\scriptsize $R$};
\node[anchor=west,font=\small] at (-7.85,-4.75) {Typed reference};
\node[ped male, ped miss] at (-3.7,-4.75) {\scriptsize $M_i$};
\node[anchor=west,font=\small] at (-3.25,-4.75) {Missing person};
\node[ped male, ped unk] at (0.9,-4.75) {};
\node[anchor=west,font=\small] at (1.35,-4.75) {Untyped relative};
\node[ped male, ped vic] at (5.5,-4.75) {\scriptsize $V_i$};
\node[anchor=west,font=\small] at (5.95,-4.75) {Post-mortem sample};
\end{tikzpicture}\caption{A small DVI problem with $s=t=4$. }\label{runex}
\end{figure}

For the sake of simplicity, we assume data on only one locus, with alleles $a,b$ and $c$. The frequencies are $p_a=0.2,p_b=0.3, p_c=0.5$, and for simplicity, we do not include any mutation model. 
We start with $F_0$. Mapping all open spots to victims results in the probability of finding the victim profiles in random unrelated individuals, as summarized in the following table. 

\begin{table}[H]
\centering
\begin{tabular}{|c|c|c|c|c|c|c|c|}
\hline
\rowcolor{gray!30}
Decision variable & $O_1$ & $O_2$ & $O_3$ & $O_4$ & Likelihood formula & Numerical value & Weight $w_{i,0}$ \\
\hline
$x_{1,0}$ & $V_1$ & $*$ & $*$ & $*$ & $2p_ap_b$ & 0.12 & $-2.12$ \\
\hline
$x_{2,0}$ & $*$ & $V_2$ &  $*$ & $*$ & $2p_ap_c$ & 0.2 & $-1.61$ \\
\hline
$x_{3,0}$ & $*$ & $*$ & $V_3$ & $*$ & $p_b^2$ & 0.09 & $-2.41$ \\
\hline
$x_{4,0}$ & $*$ & $*$ & $*$ & $V_4$ & $2p_bp_c$ & 0.3 & $-1.20$ \\
\hline
$x_{5,0}$ & $V_1$ & $V_2$ & $V_3$ & $V_4$ & $8p_a^2p_b^4p_c^2$ & 0.000648 & $-7.34$ \\
\hline
$x_{6,0}$ & $*$ & $*$ & $*$ & $*$ & $1$ & 1 & $0.00$ \\
\hline
\end{tabular}
\caption{Some assignments for $F_0$. The weights of configurations not on this list can be computed as a combination of the first four entries.}
\label{tabelF3}
\end{table}

Note that the assignment that corresponds to the decision $x_{5,0}=1$ connects all open spots with victims, thus also mapping all missing persons to untyped individuals. We call this trivial assignment the \textit{null baseline}.

Next, we calculate all feasible partial assignments for $F_1$ and $F_2$ separately. Starting with $F_1$, we note that $M_1$ must inherit $a$ from $R_1$, therefore $V_3$ or $V_4$ on $M_1$ will always result in a likelihood of 0. The weights $w_{i,1}$ of the non-zero likelihoods are presented in \Cref{tabelF1}.

\begin{table}[H]
\centering
{\begin{tabular}{|c|c|c|c|c|c|}
\hline
\rowcolor{gray!30}
Decision variable & $M_1$ & $M_2$ & Likelihood formula & Numerical value & Weight $w_{i,1}$\\
\hline
$x_{1,1}$ & $V_1$ & $V_2$ & $\frac12p_bp_c$ & 0.075 & $-2.59$ \\
\hline
$x_{2,1}$ & $V_1$ & $V_3$ & $\frac12p_b^2$ & 0.045 & $-3.10$ \\
\hline
$x_{3,1}$ & $V_1$ & $V_4$ & $\frac12p_bp_c$ & 0.075 & $-2.59$ \\
\hline
$x_{4,1}$ & $V_1$ & $*$ & $p_b$ & 0.3 & $-1.20$ \\
\hline
$x_{5,1}$ & $V_2$ & $V_1$ & $\frac12p_bp_c$ & 0.075 & $-2.59$ \\
\hline
$x_{6,1}$ & $V_2$ & $V_4$ & $\frac12p_bp_c$ & 0.075 & $-2.59$ \\
\hline
$x_{7,1}$ & $V_2$ & $*$ & $p_c$ & 0.5 & $-0.69$ \\
\hline
$x_{8,1}$ & $*$ & $V_1$ & $\frac12p_b+p_ap_b$ & 0.21 & $-1.56$ \\
\hline
$x_{9,1}$ & $*$ & $V_2$ & $\frac12p_c+p_ap_c$ & 0.35 & $-1.05$ \\
\hline
$x_{10,1}$ & $*$ & $V_3$ & $\frac12p_b^2$ & 0.045 & $-3.10$ \\
\hline
$x_{11,1}$ & $*$ & $V_4$ & $p_bp_c$ & 0.15 & $-1.90$ \\
\hline
$x_{12,1}$ & $*$ & $*$ & $1$ & 1 & $0.00$ \\
\hline
\end{tabular}}
\caption{Non-zero configurations for $F_1$ of \Cref{runex}. An $*$ indicates that a missing person is mapped to a dummy variable.}
\label{tabelF1}
\end{table}

We repeat this for $F_2$, again noting that assigning $V_1$ or $V_2$ will immediately result in likelihood 0. This yields

\begin{table}[H]
\centering
{\begin{tabular}{|c|c|c|c|c|c|}
\hline
\rowcolor{gray!30}
Decision variable & $M_3$ & $M_4$ & Likelihood formula & Numerical value & Weight\\
\hline
$x_{1,2}$ & $V_3$ & $V_4$ & $\frac12^2$ & 0.25 & $-1.39$ \\
\hline
$x_{2,2}$ & $V_3$ & $*$ & $\frac12$ & 0.5 & $-0.69$ \\
\hline
$x_{3,2}$ & $V_4$ & $*$ & $\frac12$ & 0.5 & $-0.69$ \\
\hline
$x_{4,2}$ &$*$ & $*$ & $1$ & 1 & $0.00$ \\
\hline
\end{tabular}}
\caption{Non-zero configurations for $F_2$ of \Cref{runex}. An $*$ indicates that a missing person is mapped to a dummy variable. Symmetries are not included.}
\label{tabelF2}
\end{table}

With the above tables, we can find a global optimum by combining the best partial assignments into a feasible global assignment $x$. In this case, an optimal $x$ is such that $x(M_i)=V_i$ and $x(O_i)=*_i$. In the linear program, this amounts to the decisions $x_{6,0}=x_{1,1}=x_{1,2}=1$.
    
\end{example}

\subsection{Reducing the number of local calculations}\label{pruning}
In the above example, it is possible to calculate all weights by hand. As explained above, this becomes impractical with larger families and more victims to consider. We discuss the main strategies to address this problem below. 

\subsubsection{Eliminating local assignments}
\label{eliminating_local_assignments}
As demonstrated in the example, any assignment that contains a partial assignment with zero likelihood can never be better than the null baseline of assigning all missing persons to the dummy variables. For instance, placing $V_3$ on $M_1$ will globally yield $L(x) = -\infty$ if $x$ includes this identification, since $V_3$ carries no $a$ allele and $M_1$ must inherit $a$ from reference $R_1$. This partial assignment can thus be discarded outright, without evaluating any combination of $V_3$ with other missing persons.

In order to exploit this idea more efficiently, we can make use of the assumed genetic independence between loci. Instead of immediately considering the product over all loci, we calculate the likelihoods locus by locus. As soon as a single locus already yields a zero for a partial assignment, evaluation on other loci is redundant.

When a mutation model is incorporated, which is generally the case, practically impossible partial assignments will generally not be zero but very close to zero, typically still below the null baseline. This idea still works in this case, but we have to be a bit more careful with such eliminations, since we are interested in the top $k$ best assignments. The justification lies in the observation that the null baseline is usually a very poor explanation of the data. For that reason, discarding even worse configurations does generally not affect the ultimate goal of providing a list of the most likely solutions. In the highly unlikely case that this happens, it will automatically be noted because the null baseline will then be included in the top $k$.

This principle becomes more powerful when one uses a feasible ``good guess'' instead of the null baseline. When a problem is not fully open, this is automatically required, because assigning only dummy variables is then not allowed. Feasible baselines can be easily found using current sequential algorithms \cite{vigeland2021book}.

\subsubsection{An efficient order of evaluation}
The effectiveness of finding highly unlikely partial assignments often depends on the order in which missing individuals are assigned. Preferably, this is done at an early stage, thus eliminating all downstream possibilities simultaneously. As such, the greatest gains come from assigning the most forcing individuals first. In the running example, $M_1$ is a direct child of typed reference $R_1$, meaning Mendelian inheritance strongly constrains which victim genotypes are compatible. $M_2$, by contrast, is connected to the pedigree only through $M_1$. Scanning $M_1$ first therefore eliminates far more candidates before any combination with $M_2$ is considered.


\subsubsection{Pre-filtering via inverse 
search}
So far, the improvements have been on the AM side, taking advantage of the structure of the reference families. Here, we also present an improvement for the processing of the PM data. Rather than evaluating all candidate victims against a missing person and discarding incompatible ones after the fact, we can derive the eligible pool directly from the available constraints without carrying out a full likelihood evaluation. Three sources of constraint are exploited.

\begin{itemize}
\item Typed references constrain the genotypes a missing individual can carry via Mendelian inheritance. In the running example, both parents of $M_3$ are typed ($R_2$: $b/b$, $R_3$: $b/c$), so $M_3$ can only be $(b,b)$ or $(b,c)$. Victims with any other genotype are excluded immediately. Intersecting this constraint simultaneously across multiple loci typically reduces the candidate pool from hundreds of victims to a handful. For $M_1$, with only one typed parent ($R_1$: $a/a$), only victims carrying at least one $a$ allele qualify. For $M_2$, whose only typed relative is grandparent $R_1$ separated by the unresolved slot $M_1$, no genotype constraint can be propagated without first knowing who fills $M_1$.

\item When available, the sex recorded for each missing individual in the pedigree provides an immediate filter: victims of the wrong sex are excluded at no computational cost.

\item When two missing slots stand in a parent--child relationship with each other, placing a victim in the parent slot constrains the child slot: a compatible victim must share at least one allele with the assigned parent at every typed locus. This constraint is applied dynamically as the recursion proceeds.
\end{itemize}
Together, these three filters can eliminate the vast majority of candidate victims before the likelihood evaluations even begin, leaving only a small pool of genetically plausible assignments to evaluate. As before, this idea also works when mutation is incorporated, where one discards the very unlikely local configurations with the aid of a suitable global baseline. 

Naturally, these strategies do not reveal their power in the small example above, but we will demonstrate their effectiveness in the next section. To conclude this section, we describe another, separate component of our method. It provides a way to evaluate the remaining assignments very efficiently, while also breaking symmetry along the way. In large cases this is a hugely useful tool, because it automatically removes many redundancies in the top $k$. 

\subsection{Faster evaluations}\label{sec:symmetries}
In the human case of little inbreeding and frequent application of mutation models, the standard approach for likelihood evaluations is the Elston-Stewart (ES) algorithm \cite{elston1971}. Although the current ES implementations are quite fast, we do benefit from a bespoke implementation. This involves an algorithm that returns a polynomial for any given (non-inbred) pedigree. With this polynomial, the likelihood of an assignment to a pedigree then reduces to simply plugging in zeroes and ones depending on the genetic information that we have.

This approach has two main advantages. Firstly, evaluation can be done utilizing  a \textit{graphics processing unit} (GPU) \cite{kirk2022}, offering significant advantages in case of evaluating high numbers of assignments on the same pedigree. Secondly, polynomials effectively expose symmetries such as multiple missing siblings.

Due to the high notational burden, we refer to the appendix for a complete description of the algorithm. In the remainder of this section, we only demonstrate the method on a small example.

Consider the following pedigree in \Cref{NF1} 
\tikzset{
  ped male/.style   = {draw, rectangle, minimum size=5.5mm, inner sep=0pt},
  ped female/.style = {draw, circle,    minimum size=5.5mm, inner sep=0pt},
  ped ref/.style    = {fill=black!72, text=white, thick},
  ped miss/.style   = {dashed, thick},
  ped unk/.style    = {fill=black!7},
}
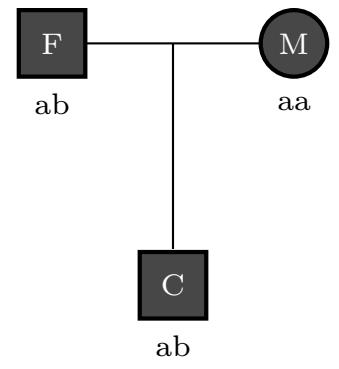
\begin{figure}
\centering
\begin{tikzpicture}[scale=0.8, every node/.style={transform shape}]

\draw (-1,1) -- (1,1);
\draw (0,1) -- (0,-0.7);
\node[ped male, ped ref] at (-1,1) {\scriptsize F}; \node at (-1,0.5) {\scriptsize ab};
\node[ped female, ped ref] at (1,1) {\scriptsize M}; \node at (1,0.5) {\scriptsize aa};
\node[ped male, ped ref] at (0,-1) {\scriptsize C}; \node at (0,-1.5) {\scriptsize ab};
\end{tikzpicture}\caption{A nuclear family with one child, with one marker shown.}\label{NF1}
\end{figure}
We work with only two alleles $a$ and $b$. This gives $K=3$ possible profiles, namely $aa, ab$ and $bb$. We further suppose that $p_a=p_b=\frac12$, so that $P(aa)=P(bb)=\frac14,P(ab)=\frac12.$

A direct calculation of the likelihood of the above assignment could run as follows. The parents $F$ and $M$ are non-inbred and thus independent, so that $P(F=ab, M=aa)=P(ab)P(aa)=\frac18$. The mother $M$ passes on $a$ to $C,$ and $C$ has a 50\% chance of inheriting $b$ from $F$, so that $P(C=ab \mid F=ab, M=aa) = \frac12.$ This yields $P(F=ab, M=aa, C=ab) = \frac1{16}.$

In our algorithmic approach, we introduce $K$ variables for each individual. In this example, we write $x_{aa}, x_{ab}, x_{bb}$ for the father. If the father is known to be $aa$ for example, we set $x_{aa}=1$ and $x_{ab}=x_{bb}=0$. If the father is unknown, we set $x_{aa}=x_{ab}=x_{bb}=1$. For the mother we write $y_{aa},y_{ab},y_{bb}$ with the same interpretation. This yields the following vector with the product of the genotype indicators and the probability of the corresponding genotypes:
\begin{equation}
    \bm{v}=\left(\frac1{16}x_{aa}y_{aa},\frac18x_{aa}y_{ab},\frac1{16}x_{aa}y_{bb},\frac18x_{ab}y_{aa},\frac14x_{ab}y_{ab},\frac18x_{ab}y_{bb},\frac1{16}x_{bb}y_{aa},\frac18x_{bb}y_{ab},\frac1{16}x_{bb}y_{bb}\right).
\end{equation}
Moving on to $C,$ we precalculate all the probabilities $P(C \mid F,M)$ and store these in a matrix (not displayed here). We write the variables $z_{aa},z_{ab},z_{bb}$ to store our knowledge about $C$. Given the information in \Cref{NF1}, we immediately conclude that $z_{aa}=z_{bb}=0, z_{ab}=1$ cf. \Cref{eliminating_local_assignments}. Thus in this example, we only need the following vector
\begin{equation}
\bm{w}=\left(0,\frac12z_{ab},z_{ab},\frac12z_{ab},\frac12z_{ab},\frac12z_{ab},z_{ab},\frac12z_{ab},0\right),
\end{equation}
where the coefficients are now equal to $P(C =ab \mid F,M)$.

The main point is that the likelihood of the entire pedigree can be obtained from the following inner product:
\begin{equation}\label{nuclinprodr}
    \langle \bm{v}, \bm{w} \rangle = \frac{z_{ab}}{16}
\left(
x_{aa}y_{ab} + x_{aa}y_{bb} + x_{ab}y_{aa} + 2x_{ab}y_{ab} + x_{ab}y_{bb} + x_{bb}y_{aa} + x_{bb}y_{ab}
\right).
\end{equation}
In case of \Cref{NF1} all variables are 0, apart from $z_{ab},x_{ab}$ and $y_{aa}.$ This yields the correct answer of $\frac1{16}$.

To demonstrate the flexibility of \eqref{nuclinprodr}, suppose for instance that the father is unknown. Then
\begin{equation}
\begin{split}
        P(M=aa, C=ab) =& ~P(C=ab, F=aa, M=aa) \\
        &+ P(C=ab, F=ab, M=aa) \\
        &+ P(C=ab, F=bb, M=aa) \\
        =&~ 0 + \frac1{16} + \frac{1}{16}=\frac1{8}.
\end{split}
\end{equation}

Indeed, when we take $x_{aa}=x_{ab}=x_{bb}=1,$ we see that
\begin{equation}
    \langle \bm{v}, \bm{w} \rangle = \frac{1}{16}
\left(
1 + 1
\right)=\frac18.
\end{equation}
This procedure translates the recursive ES-algorithm into a fixed set of (sparse) matrix operations, which are known to be very efficiently calculated on GPUs \cite{bell2008}.

\subsubsection{Symmetries}
In simple pedigrees, symmetries are usually easy to spot with the human eye. For example, it is immediately clear that in case of two missing brothers, permuting the two will not change the likelihood.

For a computer, this is much more difficult as graphs are complex objects to process. The above algebraic representations automatically yield all symmetries by identifying the symmetric variables.

As an example, if we ignore sex, it is intuitively clear that swapping the parents in a nuclear family does not change the likelihood if we only consider autosomal DNA. With our method, the algorithm can easily deduce this from \eqref{nuclinprodr}, because the polynomial is invariant under the swapping $(x_{aa}, x_{ab}, x_{bb}) \leftrightarrow (y_{aa}, y_{ab}, y_{bb})$.

A similar mechanism works when we compare different pedigrees. Suppose that in one pedigree, we are looking for a man, his father and his grandfather. In the other pedigree, we are also looking for a man and his father, but instead of the grandfather we are missing a half-brother.

\tikzset{
  ped male/.style   = {draw, rectangle, minimum size=5.5mm, inner sep=0pt},
  ped female/.style = {draw, circle,    minimum size=5.5mm, inner sep=0pt},
  ped ref/.style    = {fill=black!72, text=white, thick},
  ped miss/.style   = {dashed, thick},
  ped unk/.style    = {fill=black!7},
}
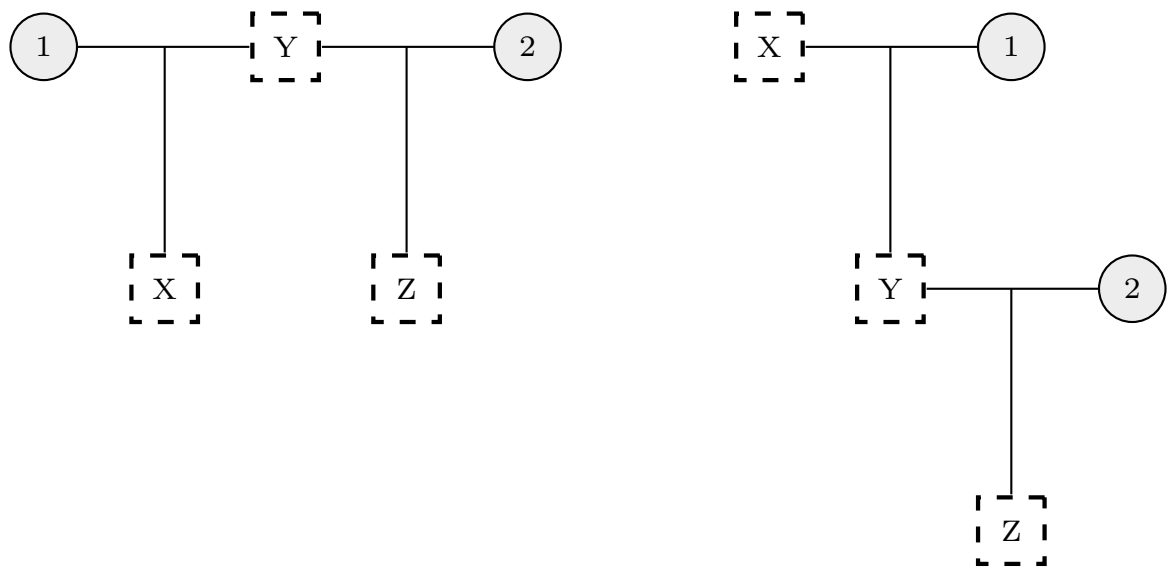
\begin{figure}[H]
\centering
\begin{tikzpicture}[scale=0.8, every node/.style={transform shape}]

\draw (-1,1) -- (0.7,1);
\draw (1.3,1) -- (3,1);
\draw (0,1) -- (0,-0.7);
\draw (2,1) -- (2, -0.7);

\node[ped female, ped unk] at (-1,1) {\scriptsize 1}; 
\node[ped male, ped miss] at (1,1) {\scriptsize Y};
\node[ped female, ped unk] at (3,1) {\scriptsize 2}; 

\node[ped male, ped miss] at (0,-1) {\scriptsize X};
\node[ped male, ped miss] at (2,-1) {\scriptsize Z};

\draw (5.3,1) -- (7,1);
\draw (6,1) -- (6,-0.7);
\draw (6.3,-1) -- (8, -1);
\draw (7,-1) -- (7,-2.7);

\node[ped male, ped miss] at (5,1) {\scriptsize X}; 
\node[ped female, ped unk] at (7,1) {\scriptsize 1};

\node[ped male, ped miss] at (6,-1) {\scriptsize Y}; 
\node[ped female, ped unk] at (8,-1) {\scriptsize 2};

\node[ped male, ped miss] at (7,-3) {\scriptsize Z};

\end{tikzpicture}\caption{Two symmetric pedigrees.}  \label{vbsymmetrie}
\end{figure}

The dependencies between $X,Y$ and $Z$ are identical in both cases. As before, this is straightforward to see from a genetic perspective, but non-obvious for a computer. However, after reducing the pedigrees to their respective inner product, the resulting polynomials will be identical.

Finally, we note that the polynomials can also be used to detect independencies within pedigrees. If the polynomial can be factorized such that all individuals are all contained in one factor, the factors represent the independent branches of the pedigree. For DVI cases however, this is generally not relevant. Large pedigrees are uncommon and recognizing and redefining separate branches of a pedigree is not difficult.
\newpage
\section{Results}\label{results}

In this section, we evaluate our method using simulated data motivated by DVI cases mentioned in the introduction. Section~\ref{sec:resexhaustive} compares our method with an exhaustive-search baseline, illustrating the combinatorial growth of the assignment space. Section~\ref{sec:12_designs} examines computational performance across twelve pedigree designs of increasing assignment space.
In section~\ref{sec:cases} we consider three case studies: a simulated passenger-flight incident involving 300 victims, a simulated mass grave involving 10,000 victims, and a quasi-empirical case concerning Argentinian children. The latter case also illustrates the interpretation of the top-$k$ highest-likelihood assignments. Finally, section~\ref{sec:ablation} presents an ablation study to assess the contribution of the algorithmic components across these designs.

For the simulations, we use 23 loci of the Norwegian STR frequency database distributed with the \textsf{forrel} package in
R \cite{dupuy2013} and we use an equal (also called uniform) mutation model \cite{egeland2016, simonsson2016}. All runs were performed on the Snellius supercomputer using a single NVIDIA A100 GPU (40\,GB) and 12 cores of an Intel Xeon Platinum 8360Y processor.

\subsection{Naive exhaustive validation on small instances}\label{sec:resexhaustive}

Consider $k$ independent families, each containing one typed grandfather and two missing grandsons, so victims and missing persons both number $2k$. See \Cref{exhaustive}.

\tikzset{
  ped male/.style   = {draw, rectangle, minimum size=5.5mm, inner sep=0pt},
  ped female/.style = {draw, circle,    minimum size=5.5mm, inner sep=0pt},
  ped ref/.style    = {fill=black!72, text=white, thick},
  ped miss/.style   = {dashed, thick},
  ped unk/.style    = {fill=black!7},
  ped vic/.style    = {fill=black!0},
}
\begin{figure}[H]
\centering
\begin{tikzpicture}[scale=0.8, every node/.style={transform shape},
  ped male/.style   = {draw, rectangle, minimum size=8.5mm, inner sep=0pt},
  ped female/.style = {draw, circle,    minimum size=8.5mm, inner sep=0pt}]

\draw (-9,2.00) -- (11,2.00);

\draw (-9,-4.00) -- (11,-4.00);

\draw (11,2.00) -- (11,-4.00);

\draw (-9,2.00) -- (-9,-4.00);

\node at (-7,1.5) {PM data};

\node[ped male, ped vic] at (-8,0.75) {\scriptsize $V_1$};
\node[ped male, ped vic] at (-8,-0.50) {\scriptsize $V_2$};
\node at (-8,-1.75) {\scriptsize $\vdots$};
\node[ped male, ped vic] at (-8,-3) {\tiny $V_{2k}$};

\node[ped male, ped vic] at (-6,0.75) {\scriptsize $*_1$};
\node[ped male, ped vic] at (-6,-0.50) {\scriptsize $*_2$};
\node at (-6,-1.75) {\scriptsize $\vdots$};
\node[ped male, ped vic] at (-6,-3) {\tiny $*_{2k}$};

\draw (-5,2.00) -- (-5,-4.00);


\node[ped male, ped miss] at (-3,0.75) {\scriptsize $O_1$};
\node[ped male, ped miss] at (-3,-0.325) {\scriptsize $O_2$};
\node at (-3,-1.5) {\scriptsize $\vdots$};
\node[ped male, ped miss] at (-3,-2.625) {\tiny $O_{2k}$};


\node at (2.5,1.5) {AM data};
\draw (0.38,0.75) -- (1.53,0.75);
\draw (0.96,0.75) -- (0.96,-0.71);
\draw (1.22,-1) -- (3.17,-1);
\draw (2.04,-1) -- (2.04,-2);
\draw (1.44,-2) -- (2.64,-2);
\draw (1.44,-2) -- (1.44,-2.33);
\draw (2.64,-2) -- (2.64,-2.33);
\node[ped male, ped ref] at (0.38,0.75) {\scriptsize $R_1$};
\node[ped female, ped unk] at (1.53,0.75) {\scriptsize $U_1$};
\node[ped male, ped unk] at (0.96,-1) {\scriptsize $U_2$};
\node[ped female, ped unk] at (3.11,-1) {\scriptsize $U_3$};
\node[ped male, ped miss] at (1.44,-2.75) {\scriptsize $M_1$};
\node[ped male, ped miss] at (2.64,-2.75) {\scriptsize $M_2$};

\draw (6.38,0.75) -- (7.53,0.75);
\draw (6.96,0.75) -- (6.96,-0.71);
\draw (7.22,-1) -- (9.17,-1);
\draw (8.04,-1) -- (8.04,-2);
\draw (7.44,-2) -- (8.64,-2);
\draw (7.44,-2) -- (7.44,-2.33);
\draw (8.64,-2) -- (8.64,-2.33);
\node[ped male, ped ref] at (6.38,0.75) {\tiny $R_{k}$};
\node[ped female, ped unk] at (7.53,0.75) {\tiny $U_{3k-2}$};
\node[ped male, ped unk] at (6.96,-1) {\scriptsize $U_{3k-1}$};
\node[ped female, ped unk] at (9.11,-1) {\tiny $U_{3k}$};
\node[ped male, ped miss] at (7.44,-2.75) {\tiny $M_{2k-1}$};
\node[ped male, ped miss] at (8.64,-2.75) {\tiny $M_{2k}$};

\node at (-3,-3.5) {$F_0$};
\node at (2,-3.5) {$F_1$};
\node at (8,-3.5) {$F_k$};

\node at (5,-1) {$\dots$};

\node[ped male, ped ref] at (-8.3,-4.75) {\scriptsize $R$};
\node[anchor=west,font=\small] at (-7.75,-4.75) {Typed reference};
\node[ped male, ped miss] at (-3.2,-4.75) {\scriptsize $M_i$};
\node[anchor=west,font=\small] at (-2.65,-4.75) {Missing person};
\node[ped male, ped unk] at (1.9,-4.75) {};
\node[anchor=west,font=\small] at (2.45,-4.75) {Untyped relative};
\node[ped male, ped vic] at (6.5,-4.75) {\scriptsize $V_i$};
\node[anchor=west,font=\small] at (7.05,-4.75) {Post-mortem sample};
\end{tikzpicture}\caption{Schematic DVI instance with \(k\) independent reference families, each containing two missing grandsons. The additional group \(F_0\) contains the open spots, while dummy victims allow missing persons to remain unidentified.}\label{exhaustive}
\end{figure}
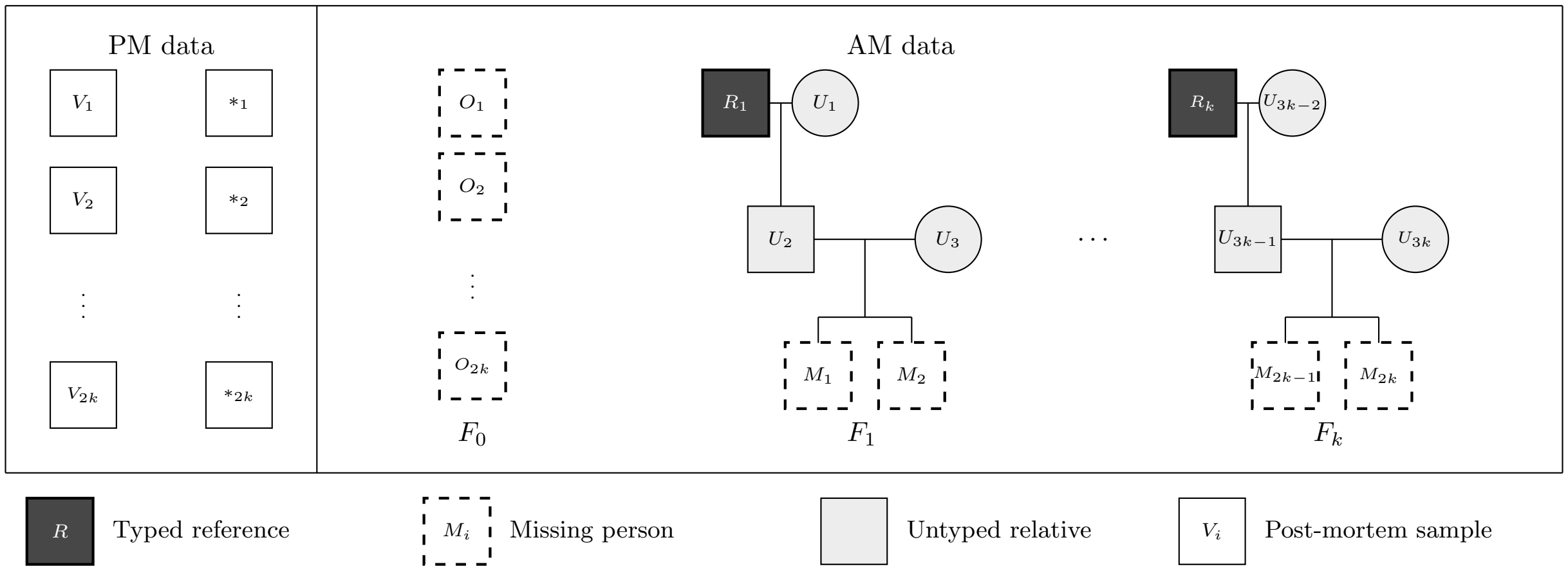

We carried this out up to $k=3$ of such reference families; since at $k \geq4$, full enumeration becomes infeasible.  
The results are displayed in \Cref{tab:exhaustive}.

\begin{table}[H]
\centering
\small
{\begin{tabular}{|l|r|r|}
\hline
\rowcolor{gray!30}
Method & Likelihood evaluations & Saved evaluations \\
\hline
\multicolumn{3}{|l|}{\cellcolor{gray!10}$k=1$: $2$ victims, $2$ missing
 ($7$ assignments)}\\
\hline
Exhaustive & $483$ & $337$ \\
Our method & $\mathbf{146}$ & $\mathbf{69{.}77\%}$ \\
\hline
\multicolumn{3}{|l|}{\cellcolor{gray!10}$k=2$: $4$ victims, $4$ missing
 ($209$ assignments)}\\
\hline
Exhaustive & $28\,842$ & $28\,395$ \\
Our method & $\mathbf{447}$ & $\mathbf{98{.}45\%}$ \\
\hline
\multicolumn{3}{|l|}{\cellcolor{gray!10}$k=3$: $6$ victims, $6$ missing
 ($13\,327$ assignments)}\\
\hline
Exhaustive & $2\,758\,689$ & $2\,757\,648$ \\
Our method & $\mathbf{1\,041}$ & $\mathbf{99{.}96\%}$ \\
\hline
\end{tabular}}
\caption{Exhaustive enumeration of distinct assignments compared with
our method. Assignments that differ only in dummy labels or open-spot
order are counted once. Savings are reported as absolute counts
and percentage reductions in likelihood evaluations relative to
exhaustive enumeration.}
\label{tab:exhaustive}
\end{table}

Across all enumerated instances, our method recovered the same maximum-likelihood solution as exhaustive search, while requiring far fewer likelihood evaluations. The relative savings increased with problem size. For \(k=3\), the exhaustive baseline requires \(2{,}758{,}689\) likelihood evaluations across \(13{,}327\) distinct assignments, compared with \(1{,}041\) evaluations for our method, corresponding to a reduction of approximately \(99.96\%\).

\subsection{Twelve designs}
\label{sec:12_designs}

The number of assignments that need to be evaluated depends on both the number of missing persons and the pedigree structure. Assignments involving parent-child relationships are genetically easier to exclude than distant relatives or full siblings. Although full siblings and parent–child pairs have the same coefficient of relationship, siblings need not share an allele at every locus, making early exclusion of incompatible candidates more difficult.

In this experiment, we define twelve pedigree structures of increasing
complexity. The pedigrees are displayed in Figure~\ref{fig:designs12}.
For each of these designs, we construct an open DVI case consisting of 10 pedigrees of that type. We apply our method and compute the average runtime after 10 repetitions. Finally, we repeat the experiment under three mutation settings: no mutation, equal mutation, and a bounded, symmetric stepwise mutation model~\cite{ohta1973, egeland2016, simonsson2016}. The average and maximum runtimes are reported in Table~\ref{tab:designs12}.

\tikzset{
  ped male/.style   = {draw, rectangle, minimum size=5.5mm, inner sep=0pt},
  ped female/.style = {draw, circle,    minimum size=5.5mm, inner sep=0pt},
  ped ref/.style    = {fill=black!72, text=white, thick},
  ped miss/.style   = {dashed, thick},
  ped unk/.style    = {fill=black!7},
}
\newcommand{\rolelabel}[3]{\node[font=\tiny] at (#1,#2) {#3};}

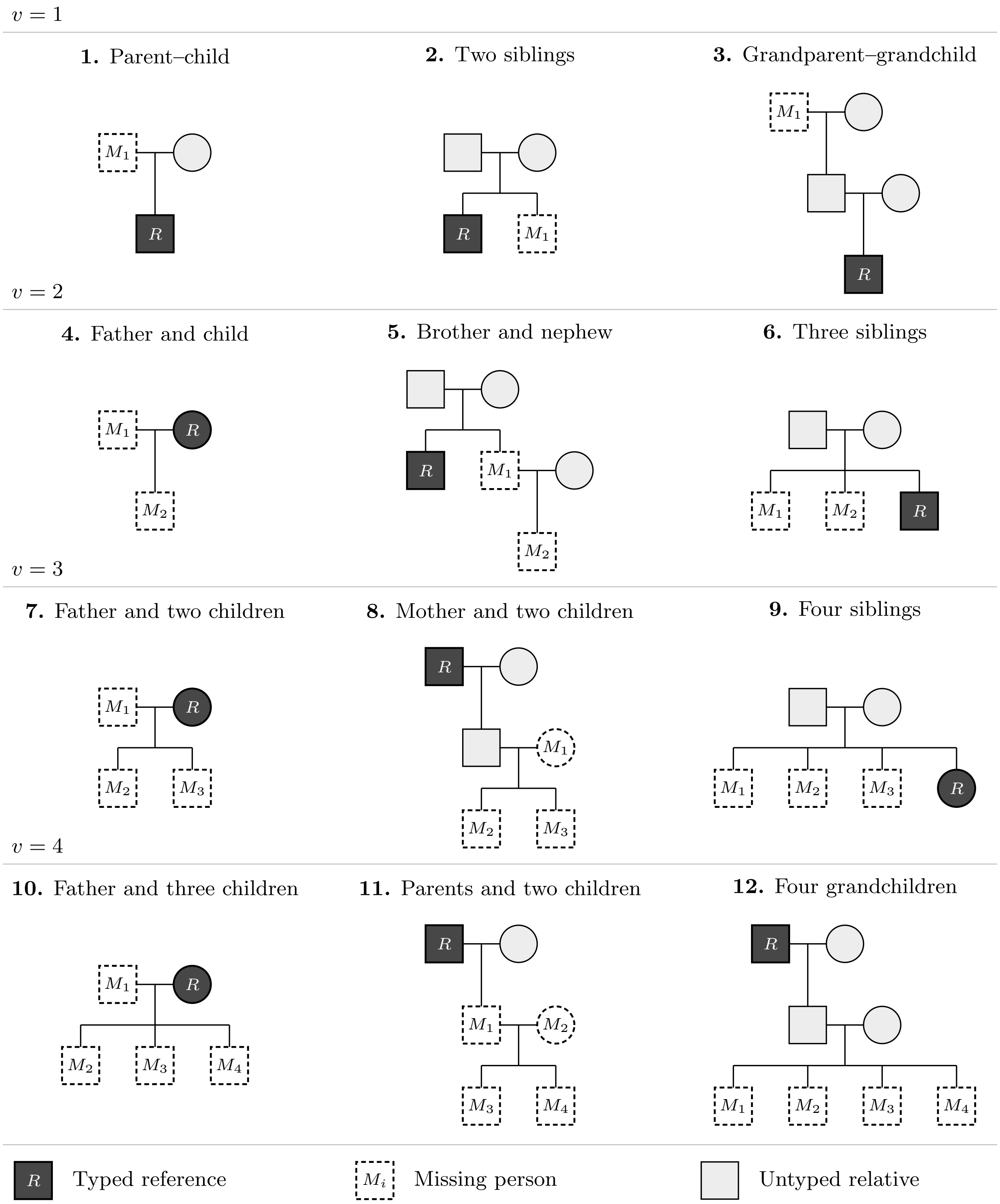
\begin{figure}
\centering
\begingroup

\newcommand{\DVIpanel}[4]{%
  \begin{scope}[shift={(#1,#2)}]
    \node[font=\small,align=center] at (0,1.48) {#3};
    #4
  \end{scope}%
}
\newcommand{\DVIunion}[3]{%
  \draw (#1.east) -- (#2.west);
  \coordinate (#3) at ($(#1.center)!0.5!(#2.center)$);
}

\resizebox{0.98\linewidth}{!}{%
\begin{tikzpicture}[
  line width=0.55pt,
  dvi person/.style={
    draw,minimum size=5.5mm,
    inner sep=0pt,outer sep=0pt,font=\scriptsize
  },
  dvi male/.style={dvi person,rectangle},
  dvi female/.style={dvi person,circle},
  dvi ref/.style={
    fill=black!72,text=white,line width=0.8pt
  },
  dvi missing/.style={
    fill=white,dash pattern=on 2pt off 1.3pt,
    line width=0.8pt
  },
  dvi untyped/.style={fill=black!7}
]

\foreach \y/\v in {0/1,-4.1/2,-8.2/3,-12.3/4}{
  \node[anchor=west,font=\small\bfseries]
    at (-2.25,\y+2.10) {$v=\v$};
  \draw[black!25,line width=0.4pt]
    (-2.25,\y+1.83) -- (12.45,\y+1.83);
}


\DVIpanel{0}{0}{\textbf{1.} Parent--child}{
  \node[dvi male,dvi missing]
    (a) at (-0.55,0.05) {$M_1$};
  \node[dvi female,dvi untyped]
    (b) at (0.55,0.05) {};
  \node[dvi male,dvi ref]
    (c) at (0,-1.15) {$R$};
  \DVIunion{a}{b}{u}
  \draw (u) -- (c.north);
}

\DVIpanel{5.1}{0}{\textbf{2.} Two siblings}{
  \node[dvi male,dvi untyped]
    (a) at (-0.55,0.05) {};
  \node[dvi female,dvi untyped]
    (b) at (0.55,0.05) {};
  \node[dvi male,dvi ref]
    (c) at (-0.55,-1.15) {$R$};
  \node[dvi male,dvi missing]
    (d) at (0.55,-1.15) {$M_1$};
  \DVIunion{a}{b}{u}
  \draw (u) -- (0,-0.55);
  \draw (-0.55,-0.55) -- (0.55,-0.55);
  \draw (-0.55,-0.55) -- (c.north);
  \draw (0.55,-0.55) -- (d.north);
}

\DVIpanel{10.2}{0}{\textbf{3.} Grandparent--grandchild}{
  \node[dvi male,dvi missing]
    (a) at (-0.825,0.65) {$M_1$};
  \node[dvi female,dvi untyped]
    (b) at (0.275,0.65) {};
  \node[dvi male,dvi untyped]
    (c) at (-0.275,-0.55) {};
  \node[dvi female,dvi untyped]
    (d) at (0.825,-0.55) {};
  \node[dvi male,dvi ref]
    (e) at (0.275,-1.75) {$R$};
  \DVIunion{a}{b}{u}
  \draw (u) -- (c.north);
  \DVIunion{c}{d}{w}
  \draw (w) -- (e.north);
}


\DVIpanel{0}{-4.1}{\textbf{4.} Father and child}{
  \node[dvi male,dvi missing]
    (a) at (-0.55,0.05) {$M_1$};
  \node[dvi female,dvi ref]
    (b) at (0.55,0.05) {$R$};
  \node[dvi male,dvi missing]
    (c) at (0,-1.15) {$M_2$};
  \DVIunion{a}{b}{u}
  \draw (u) -- (c.north);
}

\DVIpanel{5.1}{-4.1}{\textbf{5.} Brother and nephew}{
  \node[dvi male,dvi untyped]
    (a) at (-1.10,0.65) {};
  \node[dvi female,dvi untyped]
    (b) at (0,0.65) {};
  \node[dvi male,dvi ref]
    (c) at (-1.10,-0.55) {$R$};
  \node[dvi male,dvi missing]
    (d) at (0,-0.55) {$M_1$};
  \node[dvi female,dvi untyped]
    (e) at (1.10,-0.55) {};
  \node[dvi male,dvi missing]
    (f) at (0.55,-1.75) {$M_2$};
  \DVIunion{a}{b}{u}
  \draw (u) -- (-0.55,0.05);
  \draw (-1.10,0.05) -- (0,0.05);
  \draw (-1.10,0.05) -- (c.north);
  \draw (0,0.05) -- (d.north);
  \DVIunion{d}{e}{w}
  \draw (w) -- (f.north);
}

\DVIpanel{10.2}{-4.1}{\textbf{6.} Three siblings}{
  \node[dvi male,dvi untyped]
    (a) at (-0.55,0.05) {};
  \node[dvi female,dvi untyped]
    (b) at (0.55,0.05) {};
  \node[dvi male,dvi missing]
    (c) at (-1.10,-1.15) {$M_1$};
  \node[dvi male,dvi missing]
    (d) at (0,-1.15) {$M_2$};
  \node[dvi male,dvi ref]
    (e) at (1.10,-1.15) {$R$};
  \DVIunion{a}{b}{u}
  \draw (u) -- (0,-0.55);
  \draw (-1.10,-0.55) -- (1.10,-0.55);
  \draw (-1.10,-0.55) -- (c.north);
  \draw (0,-0.55) -- (d.north);
  \draw (1.10,-0.55) -- (e.north);
}


\DVIpanel{0}{-8.2}{\textbf{7.} Father and two children}{
  \node[dvi male,dvi missing]
    (a) at (-0.55,0.05) {$M_1$};
  \node[dvi female,dvi ref]
    (b) at (0.55,0.05) {$R$};
  \node[dvi male,dvi missing]
    (c) at (-0.55,-1.15) {$M_2$};
  \node[dvi male,dvi missing]
    (d) at (0.55,-1.15) {$M_3$};
  \DVIunion{a}{b}{u}
  \draw (u) -- (0,-0.55);
  \draw (-0.55,-0.55) -- (0.55,-0.55);
  \draw (-0.55,-0.55) -- (c.north);
  \draw (0.55,-0.55) -- (d.north);
}

\DVIpanel{5.1}{-8.2}{\textbf{8.} Mother and two children}{
  \node[dvi male,dvi ref]
    (a) at (-0.825,0.65) {$R$};
  \node[dvi female,dvi untyped]
    (b) at (0.275,0.65) {};
  \node[dvi male,dvi untyped]
    (c) at (-0.275,-0.55) {};
  \node[dvi female,dvi missing]
    (d) at (0.825,-0.55) {$M_1$};
  \node[dvi male,dvi missing]
    (e) at (-0.275,-1.75) {$M_2$};
  \node[dvi male,dvi missing]
    (f) at (0.825,-1.75) {$M_3$};
  \DVIunion{a}{b}{u}
  \draw (u) -- (c.north);
  \DVIunion{c}{d}{w}
  \draw (w) -- (0.275,-1.15);
  \draw (-0.275,-1.15) -- (0.825,-1.15);
  \draw (-0.275,-1.15) -- (e.north);
  \draw (0.825,-1.15) -- (f.north);
}

\DVIpanel{10.2}{-8.2}{\textbf{9.} Four siblings}{
  \node[dvi male,dvi untyped]
    (a) at (-0.55,0.05) {};
  \node[dvi female,dvi untyped]
    (b) at (0.55,0.05) {};
  \node[dvi male,dvi missing]
    (c) at (-1.65,-1.15) {$M_1$};
  \node[dvi male,dvi missing]
    (d) at (-0.55,-1.15) {$M_2$};
  \node[dvi male,dvi missing]
    (e) at (0.55,-1.15) {$M_3$};
  \node[dvi female,dvi ref]
    (f) at (1.65,-1.15) {$R$};
  \DVIunion{a}{b}{u}
  \draw (u) -- (0,-0.55);
  \draw (-1.65,-0.55) -- (1.65,-0.55);
  \draw (-1.65,-0.55) -- (c.north);
  \draw (-0.55,-0.55) -- (d.north);
  \draw (0.55,-0.55) -- (e.north);
  \draw (1.65,-0.55) -- (f.north);
}


\DVIpanel{0}{-12.3}{\textbf{10.} Father and three children}{
  \node[dvi male,dvi missing]
    (a) at (-0.55,0.05) {$M_1$};
  \node[dvi female,dvi ref]
    (b) at (0.55,0.05) {$R$};
  \node[dvi male,dvi missing]
    (c) at (-1.10,-1.15) {$M_2$};
  \node[dvi male,dvi missing]
    (d) at (0,-1.15) {$M_3$};
  \node[dvi male,dvi missing]
    (e) at (1.10,-1.15) {$M_4$};
  \DVIunion{a}{b}{u}
  \draw (u) -- (0,-0.55);
  \draw (-1.10,-0.55) -- (1.10,-0.55);
  \draw (-1.10,-0.55) -- (c.north);
  \draw (0,-0.55) -- (d.north);
  \draw (1.10,-0.55) -- (e.north);
}

\DVIpanel{5.1}{-12.3}{\textbf{11.} Parents and two children}{
  \node[dvi male,dvi ref]
    (a) at (-0.825,0.65) {$R$};
  \node[dvi female,dvi untyped]
    (b) at (0.275,0.65) {};
  \node[dvi male,dvi missing]
    (c) at (-0.275,-0.55) {$M_1$};
  \node[dvi female,dvi missing]
    (d) at (0.825,-0.55) {$M_2$};
  \node[dvi male,dvi missing]
    (e) at (-0.275,-1.75) {$M_3$};
  \node[dvi male,dvi missing]
    (f) at (0.825,-1.75) {$M_4$};
  \DVIunion{a}{b}{u}
  \draw (u) -- (c.north);
  \DVIunion{c}{d}{w}
  \draw (w) -- (0.275,-1.15);
  \draw (-0.275,-1.15) -- (0.825,-1.15);
  \draw (-0.275,-1.15) -- (e.north);
  \draw (0.825,-1.15) -- (f.north);
}

\DVIpanel{10.2}{-12.3}{\textbf{12.} Four grandchildren}{
  \node[dvi male,dvi ref]
    (a) at (-1.10,0.65) {$R$};
  \node[dvi female,dvi untyped]
    (b) at (0,0.65) {};
  \node[dvi male,dvi untyped]
    (c) at (-0.55,-0.55) {};
  \node[dvi female,dvi untyped]
    (d) at (0.55,-0.55) {};
  \node[dvi male,dvi missing]
    (e) at (-1.65,-1.75) {$M_1$};
  \node[dvi male,dvi missing]
    (f) at (-0.55,-1.75) {$M_2$};
  \node[dvi male,dvi missing]
    (g) at (0.55,-1.75) {$M_3$};
  \node[dvi male,dvi missing]
    (h) at (1.65,-1.75) {$M_4$};
  \DVIunion{a}{b}{u}
  \draw (u) -- (c.north);
  \DVIunion{c}{d}{w}
  \draw (w) -- (0,-1.15);
  \draw (-1.65,-1.15) -- (1.65,-1.15);
  \draw (-1.65,-1.15) -- (e.north);
  \draw (-0.55,-1.15) -- (f.north);
  \draw (0.55,-1.15) -- (g.north);
  \draw (1.65,-1.15) -- (h.north);
}

\draw[black!25,line width=0.4pt]
  (-2.25,-14.62) -- (12.45,-14.62);

\node[dvi male,dvi ref]
  at (-1.80,-15.15) {$R$};
\node[anchor=west,font=\small]
  at (-1.35,-15.15) {Typed reference};

\node[dvi male,dvi missing]
  at (3.25,-15.15) {$M_i$};
\node[anchor=west,font=\small]
  at (3.70,-15.15) {Missing person};

\node[dvi male,dvi untyped]
  at (8.35,-15.15) {};
\node[anchor=west,font=\small]
  at (8.80,-15.15) {Untyped relative};

\end{tikzpicture}%
}
\endgroup
\par\vspace{3mm}

\caption{The twelve pedigree designs, grouped by the number of missing
persons $v$. 
}
\label{fig:designs12}
\end{figure}

\begin{table}[H]
\centering
\small
\begin{tabular}{|c|c|l|r|r|r|r|r|r|}
\hline
\rowcolor{gray!30}
 & & & \multicolumn{2}{c|}{No mutation} & \multicolumn{2}{c|}{Equal } & \multicolumn{2}{c|}{Stepwise} \\
\rowcolor{gray!30}
Design & $v$ & Strongest relationship & mean & worst & mean & worst & mean & worst \\
\hline
 1 & 1 & Parent-Child    & 3.5\,s & 6.1\,s & 1.2\,s & 3.5\,s & 1.0\,s & 2.8\,s \\
 2 & 1 & Sibling   & 4.6\,s & 7.4\,s & 2.3\,s & 9.0\,s & 2.9\,s & 4.0\,s \\
 3 & 1 & Grandparent-Grandchild & 4.0\,s & 5.3\,s & 1.3\,s & 1.8\,s & 4.5\,s & 5.1\,s \\
\hline
 4 & 2 & Parent-Child     & 4.6\,s & 6.3\,s & 1.7\,s & 2.5\,s & 5.9\,s & 7.6\,s \\
 5 & 2 & Sibling/Parent-Child   & 6.5\,s & 7.9\,s & 6.1\,s & 8.2\,s & 26.4\,s & 32.0\,s \\
 6 & 2 & Sibling   & 5.2\,s & 6.5\,s & 5.1\,s & 6.9\,s & 51.5\,s & 86.8\,s \\
\hline
 7 & 3 & Parent-Child     & 6.1\,s & 7.9\,s & 3.2\,s & 4.2\,s & 15.0\,s & 17.0\,s \\
 8 & 3 & (Grand)parent-(Grand)child & 12.5\,s & 14.5\,s & 18.1\,s & 31.9\,s & 112.6\,s & 175.7\,s \\
 9 & 3 & Sibling   & 8.2\,s & 10.4\,s & 49.6\,s & 60.1\,s & 161.0\,s & 217.6\,s \\
\hline
10 & 4 & Parent-Child     & 8.3\,s & 10.3\,s & 6.0\,s & 7.8\,s & 32.1\,s & 35.6\,s \\
11 & 4 & Parent-Child & 18.1\,s & 32.1\,s & 63.4\,s & 79.1\,s & 339.6\,s & 481.2\,s \\
12 & 4 & Grandparent-Grandchild/Siblings & 44.5\,s & 66.2\,s & 292.2\,s & 423.2\,s & 999.5\,s & 1512.0\,s \\
\hline
\end{tabular}
\caption{\textbf{Solve time for the twelve designs.} Mean and worst-case
total solve time in seconds over
ten replicates. Equal: mutation rate $0.001$; Stepwise: $p_1=0.0012$,
$p_2=0.0001$. 
}
\label{tab:designs12}
\end{table}

Not surprisingly, the runtime increases when the family relationships become weaker. Notably, the more flexible equal mutation model is often faster than the stepwise mutation model. This is due to the fact that allele lumping 
is often more effective in the case of an equal mutation model. Without a mutation model, certain additional optimizations are possible which are more effective for computationally more demanding cases, however, they do come at a slight cost for computationally easier cases.

\subsection{Case studies}\label{sec:cases}

In this section we discuss three case studies.

\subsubsection{A passenger flight}
We demonstrate our method on simulated data that represents a plane crash involving 300 individuals (the Malaysia Airlines Flight 17 involved 298 passengers and crew). To the best of our knowledge, no methods exist to jointly solve a DVI problem of this scale.

We simulate a DVI case with 100 pedigrees. Each pedigree has at most 4 missing members. Victim relationships are generally close: after a randomly chosen first victim, each further victim is drawn with probability 0.5 from the immediate relatives (parents, children or full siblings) of those already chosen, and otherwise uniformly over the remainder of the pedigree. Sampling victims in this manner resembles the characteristics of a passenger flight: distant family relationships of fellow travelers are uncommon. 
Our method solved this instance in 4.5 minutes.

\subsubsection{Mass grave}
In this example, we consider an instance resembling a mass grave. We simulate this similarly to the plane crash. However, we consider a maximum of 3 victims per pedigree. The instance comprises 3500 pedigrees with 5132 missing persons, plus 3000 victims without a reference family and 3000 missing persons in reference pedigrees without a victim. The problem is  open on both sides.

Our method solves this instance in 29 hours. The computational cost is highly skewed: the median pedigree takes 7.7 seconds, whereas the five most expensive pedigrees account for 46\% of the compute needed and the four heaviest take over 26 hours each. These are pedigrees in which three members are missing and no typed parent or child constrains the search, so that the number of candidate assignments explodes. In practice, these pedigrees could be handled separately or using different methods to significantly shorten the computation time.

Note that the parameters of this instance are very similar to large disasters such as earthquakes and tsunamis.

\subsubsection{Argentinian children}
We constructed a DVI case based on the 25 pedigrees in Supplementary Table S1 of \cite{kling2017}, with one missing person per family. Using 23 STR loci and Norwegian allele frequencies, we simulated AM reference genotypes and 25 PM profiles. The families were analysed jointly while assuming genetic independence between pedigrees. The highest-ranked assignment had a substantially larger log-likelihood $-9509$ than the second-ranked assignment $-9549$, indicating strong relative support for the first assignment with a likelihood ratio of $2.61 \times 10^{17}$.

\subsection{Ablation study}
\label{sec:ablation}
To assess which components drive the framework’s computational gains, we performed a leave-one-component-out ablation on the twelve pedigree designs in \Cref{fig:designs12}. For each instance, we changed one component or execution setting at a time, keeping the simulated data and all remaining settings fixed. Table \ref{tab:full-baseline-runtimes} reports runtimes relative to the complete method for variants targeting search reductions, precomputation, GPU acceleration, and worker parallelism. These comparisons quantify each component’s incremental contribution within the complete pipeline and how that contribution varies across pedigree structures.

\newcommand{\rt}[3]{%
  \cellcolor{black!#1}%
  \ifnum#1>50\color{white}\fi%
  #2\,(\ensuremath{#3\times})%
}
\newcommand{\rtbold}[3]{%
  \cellcolor{black!#1}%
  \ifnum#1>50\color{white}\fi%
  \bfseries\boldmath%
  #2\,(\ensuremath{#3\times})%
}
\begin{table}
\centering
\scriptsize
\setlength{\tabcolsep}{3pt}
\renewcommand{\arraystretch}{1.30}
\resizebox{\linewidth}{!}{%
\begin{tabular}{lrr*{9}{r}}
\toprule
Design & $v$ & $V$ & \textbf{full} & \shortstack{No\\Lumping} & \shortstack{No PC\\Filter} & \shortstack{No\\Precompute} & \shortstack{No P1\\filter} & \shortstack{No Obj\\Filter} & \shortstack{No\\Symmetry} & CPU & No parallelism \\
\midrule
1 & 1 & 240 & \rt{0}{14}{1.00} & \rt{23}{43}{3.10} & \rt{0}{10}{0.74} & \rt{0}{10}{0.69} & \rt{0}{8}{0.60} & \rt{0}{11}{0.80} & \rt{5}{20}{1.45} & \rt{0}{9}{0.67} & \rt{0}{14}{0.99} \\
2 & 1 & 240 & \rt{0}{16}{1.00} & \rt{24}{51}{3.15} & \rt{0}{12}{0.71} & \rt{0}{13}{0.83} & \rt{0}{11}{0.70} & \rt{3}{21}{1.27} & \rt{10}{31}{1.89} & \rt{0}{12}{0.72} & \rt{22}{48}{2.95} \\
3 & 1 & 240 & \rt{0}{27}{1.00} & \rtbold{49}{149}{5.43} & \rt{0}{27}{1.00} & \rt{0}{28}{1.01} & \rt{0}{27}{0.97} & \rt{0}{25}{0.92} & \rt{9}{50}{1.84} & \rt{0}{27}{0.97} & \rtbold{72}{204}{7.45} \\
4 & 2 & 240 & \rt{0}{66}{1.00} & \rt{13}{146}{2.21} & \rt{0}{68}{1.03} & \rt{0}{56}{0.84} & \rt{0}{54}{0.81} & \rt{4}{87}{1.32} & \rt{4}{87}{1.32} & \rt{18}{175}{2.66} & \rt{19}{181}{2.75} \\
5 & 2 & 240 & \rt{0}{105}{1.00} & \rt{23}{319}{3.05} & \rt{1}{115}{1.10} & \rt{0}{107}{1.02} & \rt{0}{108}{1.03} & \rt{4}{146}{1.40} & \rt{7}{174}{1.66} & \rt{25}{339}{3.24} & \rt{28}{371}{3.54} \\
6 & 2 & 240 & \rt{0}{72}{1.00} & \rt{13}{154}{2.13} & \rt{0}{74}{1.02} & \rt{0}{73}{1.02} & \rt{1}{80}{1.11} & \rt{10}{137}{1.90} & \rt{13}{156}{2.16} & \rt{39}{328}{4.55} & \rt{36}{302}{4.20} \\
7 & 3 & 240 & \rt{0}{75}{1.00} & \rt{25}{245}{3.27} & \rt{13}{162}{2.17} & \rt{0}{67}{0.89} & \rt{0}{66}{0.88} & \rt{43}{366}{4.89} & \rtbold{100}{748}{10.01} & \rt{36}{320}{4.28} & \rt{29}{270}{3.61} \\
8 & 3 & 240 & \rt{0}{127}{1.00} & \rt{10}{245}{1.94} & \rt{1}{143}{1.13} & \rt{0}{109}{0.86} & \rt{0}{119}{0.94} & \rt{21}{362}{2.86} & \rtbold{100}{1265}{10.00} & \rt{37}{552}{4.36} & \rt{25}{413}{3.26} \\
9 & 3 & 240 & \rt{0}{208}{1.00} & \rt{66}{1394}{6.70} & \rt{0}{207}{0.99} & \rt{10}{395}{1.90} & \rt{0}{217}{1.04} & \rtbold{45}{1053}{5.06} & \rtbold{100}{2078}{10.00} & \rtbold{57}{1268}{6.10} & \rt{13}{455}{2.19} \\
10 & 4 & 240 & \rt{0}{466}{1.00} & \rt{32}{1789}{3.84} & \rtbold{100}{4658}{10.00} & \rt{0}{346}{0.74} & \rt{0}{447}{0.96} & \rtbold{100}{4658}{10.00} & \rtbold{100}{4658}{10.00} & \rt{40}{2139}{4.59} & \rt{39}{2094}{4.50} \\
11 & 4 & 240 & \rt{0}{361}{1.00} & \rt{25}{1169}{3.24} & \rtbold{55}{2136}{5.91} & \rt{0}{241}{0.67} & \rt{0}{345}{0.96} & \rtbold{66}{2513}{6.96} & \rtbold{100}{3613}{10.00} & \rtbold{53}{2082}{5.76} & \rt{42}{1733}{4.80} \\
12 & 4 & 240 & \rt{0}{5312}{1.00} & \rtbold{60}{34194}{6.44} & \rt{0}{5504}{1.04} & \rt{8}{9153}{1.72} & \rt{0}{5406}{1.02} & \rtbold{100}{53116}{10.00} & \rtbold{100}{53116}{10.00} & \rtbold{83}{44800}{8.43} & \rt{15}{12651}{2.38} \\
\bottomrule
\end{tabular}%
}
\caption{Runtimes in seconds (rounded to integers), with ratios to the \texttt{full} run in parentheses. $v$ denotes the number of victims per design and $V$ the total number of victims in the DVI case. The uniform mutation model is used. Ratios and shading use unrounded data. \texttt{full} and faster runs are white; grayscale increases linearly to black at 10$\times$ baseline, which was set as a time-out condition.}
\label{tab:full-baseline-runtimes}
\end{table}

Table \ref{tab:full-baseline-runtimes} shows that the effectiveness of the different strategies depends on pedigree structure. For example in design 11, assigning the one missing parent immediately constrains their children and hence the other parent, whereas in design 12 the missing siblings' parents are untyped, leaving no direct parent–child pairs to check. Removing parent–child filtering therefore slows design 11 by 5.91 but barely affects design 12 (1.04). Symmetry reduction, cache sharing, and objective filtering matter more as the search space grows. Allele lumping affects the time needed for the likelihood evaluations and thus impacts all designs.

The compute ablations show how the bottleneck shifts with problem size. On the smallest workloads (designs 1–3), CPU-only runs are competitive because GPU launch and transfer overhead offsets the speedup. On designs 9 and 12, however, dropping to CPU slows execution by  6.10 and 8.43 respectively, whereas dropping to a single worker costs only 2.19 and 2.38 respectively. Because early pruning is weaker in these designs, more likelihood evaluations reach the GPU. Once the shared device is fully used, adding workers helps very little. Similarly, precomputation pays off in designs 9 and 12 where the setup cost is amortized due to the large number of evaluations and the relatively expensive likelihood evaluations. The effectiveness of the P1 filter is strongly dependent on the mutation model used. For equal mutation, it does not provide any speedup as no branch can be cut early.

\section{Discussion}\label{discussion}
Although this paper is completely built on the presumption that DVI cases should be analysed and reported jointly, this is not always possible in practice. In large cases such as the Malaysia Airlines Flight 17, the circumstances often naturally lead to a sequential approach. Some victims are immediately found and readily identified, while other samples are found much later. Another factor is the quality of the samples: sometimes new techniques have to be developed in order to retrieve suitable DNA profiles.

In the meantime, the (usually strong) evidence at hand is used to identify as many victims as possible. This process is sequential in nature. However, because our method runs so fast, it can be used multiple times during the identification process. One can add as many dummy variables as needed to jointly run the analysis for the recovered victim samples and rerun the algorithm as more victims come in. In this way, new identifications can be made while simultaneously updating or validating earlier conclusions.

Even when the reporting process is carried out sequentially, it is possible to use our method as an ultimate and global check once all victims are found.

Several methodological caveats apply more broadly to kinship-based identification and fall outside the scope of this paper. Assumptions underlying standard kinship calculations, including Hardy--Weinberg equilibrium and marker independence, have been extensively studied elsewhere \cite{balding2015}, as have the complications introduced by degraded or low-quality DNA samples \cite{tvedebrink2012}.

A related concern specific to large-scale DVI operations is the choice of allele frequency database. Because victims and missing persons may originate from different populations, which naturally occur on large commercial planes, no single database may be clearly appropriate, and results can be sensitive to this choice \cite{fortier2020}. When this is a concern, we recommend running the analysis under several candidate databases: agreement across these sensitivity checks provides added confidence that the conclusions do not meaningfully depend on a particular choice of reference population.


Although the method is aimed at DVI cases, there may be more applications for which the computational framework can be useful. For example, instances of \textit{familial searching} \cite{slooten2014} can be reformulated as DVI problems.

\section{Conclusion}
We have developed a framework that makes joint DNA-based DVI feasible in large cases by providing methods to reduce the computational burden that was pointed out by \cite{vigeland2021}. This is achieved by combining local pedigree likelihood calculations with global integer linear optimization. Genetic filtering, likelihood-based pruning, symmetry reduction, and GPU evaluation have shown to substantially reduce the runtime. We have demonstrated that typical instances with up to 10.000 victims can be readily solved. Its practical performance depends strongly on pedigree structure and the availability of informative reference relatives.

\appendix

\section{Proof of complexity}\label{npbewijs}
In this proof, we do not incorporate mutations and assume that $s=t, s'=t'=0$. This is a special case of the general problem and it suffices to prove hardness for such cases. We consider the following decision variant.

\begin{itemize}
\item \emph{Instance:} A set of typed PM samples $\mathcal{V}=\{V_1,\dots,V_n\}$,
a set of missing individuals $\mathcal{M}=\{M_1,\dots,M_n\}$ with an
associated pedigree and AM genetic data.
\item \emph{Question:} Does there exist an assignment
$x:\mathcal{M}\to\mathcal{V}$ such that $L(x)>-\infty$?
\end{itemize}

We show that this problem is NP-complete by proving (i)~membership in NP
and (ii)~NP-hardness via a polynomial-time reduction from the Hamiltonian
Path problem \cite{karp1972}.

\subsection{Membership in NP}

To establish $\mathrm{DVI}\in\mathrm{NP}$, consider a certificate consisting
of an assignment $x$.

Given $x$, the likelihood $L(x)$ factorizes over the pedigrees induced by
the assignment. For each pedigree, the likelihood is computed using Elston--Stewart. This
computation is polynomial in the number of loci and linear in the number
of pedigree members (\cite{vigeland2021book}).

Since the total number of pedigrees is at most of order $|\mathcal{V}|$, the overall log-likelihood $L(x)$ can be evaluated in time polynomial in the input size. Checking if $L(x)>-\infty$ is trivial. Hence a proposed assignment can be verified in polynomial time, and therefore $\mathrm{DVI}\in\mathrm{NP}$.

\subsection{NP-hardness}

To prove NP-hardness, we provide a Karp reduction \cite{karp1972} from the
Hamiltonian Path (HP) problem, which reads as follows:
\begin{itemize}
\item \emph{Instance:} An undirected graph $G=(W,E)$ with $|W|=n$.
\item \emph{Question:} Does there exist a permutation
$(v_{\pi(1)},\dots,v_{\pi(n)})$ such that
$(v_{\pi(i)},v_{\pi(i+1)})\in E$ for all $1\le i<n$?
\end{itemize}

Given a graph $G=(W,E)$ with $|W|=n$, we construct in polynomial time a
DVI instance $f(G)$ as follows.

\subsubsection{Pedigree structure.}
We construct a linear pedigree consisting of $n$ missing individuals
$M_1,\dots,M_n$, where for each $i=1,\dots,n-1$, the individual $M_{i+1}$
is the child of $M_i$ and a typed reference spouse $R_i$. The spouses
$R_i$ are unrelated founders. This is a single reference family $F_1$
with $I_1=\{M_1,\dots,M_n\}$ and $d_1=n$.

\subsubsection{PM samples and missing persons.}
We create $n$ typed PM samples $V_1,\dots,V_n$, establishing a one-to-one
correspondence between vertices $W=\{v_1,\dots,v_n\}$ and the missing
individuals. An assignment therefore induces an ordering of the vertices
along the linear pedigree.

\subsubsection{Genetic construction.}
For every unordered pair of distinct vertices $\{v_j,v_k\}$ with
$\{v_j,v_k\}\notin E$ we introduce a locus $\ell_{jk}$ with two alleles
$\{a,b\}$, with $p_a=p_b=0.5$.

At locus $\ell_{jk}$:
\begin{itemize}
\item Victim $V_j$ is assigned genotype $aa$.
\item Victim $V_k$ is assigned genotype $bb$.
\item All other victims and all reference spouses $R_i$ are assigned
genotype $ab$.
\end{itemize}

We show that $G$ contains a Hamiltonian path if and only if the constructed
DVI instance admits an assignment with strictly positive likelihood.

\paragraph{($\Rightarrow$)}
Suppose $G$ contains a Hamiltonian path $v_{\pi(1)}\to\dots\to v_{\pi(n)}$.
Assign the corresponding PM samples to the missing individuals
$M_1,\dots,M_n$ in this order, i.e.\ $x(M_i)=V_{\pi(i)}$.

All members are typed once $x$ is fixed, so the likelihood at a locus is
the product of the founder probabilities of $M_1,R_1,\dots,R_{n-1}$ and
the transmission probabilities $P(g_{M_{i+1}}\mid g_{M_i},g_{R_i})$. The
founder probabilities are positive. For every adjacent pair
$(v_{\pi(i)},v_{\pi(i+1)})$, an edge exists in $E$, so no locus $\ell_{jk}$
corresponding to a non-edge forces opposite homozygotes on a parent--child
pair. Consequently, the overall likelihood is strictly positive.

\paragraph{($\Leftarrow$)}
Conversely, suppose there exists an assignment $x$ with $L(x)>-\infty$. This implies that no Mendelian incompatibility occurs at any locus. Because $R_i$ has genotype $ab$ and can transmit either allele, a trio $(M_i,R_i,M_{i+1})$ is inconsistent exactly when parent $M_i$ and child $M_{i+1}$ are opposite homozygotes. At locus $\ell_{jk}$ this happens if and only if $V_j$ and $V_k$ are placed on consecutive members of the chain, in either order. Therefore every consecutive pair in the induced ordering must correspond to an edge in $E$, implying that the assignment defines a Hamiltonian path in $G$.

The construction introduces at most $O(n^2)$ loci and $O(n)$ individuals. Every single locus specification requires constant time, yielding an overall construction time of $O(n^3)$. Therefore the reduction is computable in polynomial time.

The DVI decision problem belongs to NP and is NP-hard via polynomial-time reduction from Hamiltonian Path. Since Hamiltonian Path is NP-complete, it follows that the DVI decision problem is NP-complete. Consequently, the associated likelihood maximization problem is NP-hard.

\subsection{Example construction}
Below, we provide an example of the construction defined above.

\tikzset{
  hp vertex/.style={draw,circle,thick,minimum size=8mm,inner sep=0pt},
  hp missing/.style={draw,rectangle,dashed,thick,
                    minimum width=10mm,minimum height=7mm,inner sep=0pt},
  hp ref/.style={draw,circle,fill=black!72,text=white,thick,
                 minimum size=8mm,inner sep=0pt},
}

\begin{figure}[H]
\centering
\resizebox{.7\linewidth}{!}{%
\begin{tikzpicture}[font=\small,line cap=round]
  \node[font=\bfseries] at (2.9,7.1) {(a) Hamiltonian-path graph};
  \node[hp vertex] (v1) at (0.9,4.25) {$v_1$};
  \node[hp vertex] (v2) at (2.1,6.05) {$v_2$};
  \node[hp vertex] (v3) at (3.35,4.25) {$v_3$};
  \node[hp vertex] (v4) at (5.1,4.25) {$v_4$};
  \foreach \a/\b in {v1/v2,v1/v3,v2/v3,v3/v4}
    \draw[thick] (\a)--(\b);
  \node[font=\footnotesize] at (2.9,2.75)
    {Hamiltonian path: $v_1\!\to\!v_2\!\to\!v_3\!\to\!v_4$};
  \node[font=\footnotesize] at (2.9,2.15)
    {Nonedges: $\{v_1,v_4\}$ and $\{v_2,v_4\}$};

  \node[font=\bfseries] at (10.5,7.1) {(b) Corresponding DVI pedigree};
  \node[hp missing] (U1) at (8.45,6.05) {$U_1$};
  \node[hp ref]     (R1) at (10.35,6.05) {$R_1$};
  \node[hp missing] (U2) at (9.40,4.60) {$U_2$};
  \node[hp ref]     (R2) at (11.30,4.60) {$R_2$};
  \node[hp missing] (U3) at (10.35,3.15) {$U_3$};
  \node[hp ref]     (R3) at (12.25,3.15) {$R_3$};
  \node[hp missing] (U4) at (11.30,1.70) {$U_4$};
  \draw (U1.east)--(R1.west);
  \draw (9.40,6.05)--(U2.north);
  \draw (U2.east)--(R2.west);
  \draw (10.35,4.60)--(U3.north);
  \draw (U3.east)--(R3.west);
  \draw (11.30,3.15)--(U4.north);

\node[font=\bfseries] at (7.2,0.62)
  {(c) PM and reference genotypes};

\begin{scope}[shift={(4.17,0.30)}]

  \fill[black!12] (0,0)
    rectangle (6.06,-0.54);

  \draw[semithick] (0,0)
    rectangle (6.06,-3.24);

  \draw (1.9,0)--(1.9,-3.24);
  \foreach \c in {1,...,3}{
    \draw ({1.9+\c*1.04},0)
      --({1.9+\c*1.04},-3.24);
  }

  \foreach \r in {1,...,5}{
    \draw (0,{-0.54*\r})
      --(6.06,{-0.54*\r});
  }

  \node[font=\scriptsize] at (0.95,-0.27)
    {Profile};

  \foreach \c/\jk in {
    1/14,2/24,3/41,4/42
  }{
    \node[font=\scriptsize]
      at ({1.9+(\c-.5)*1.04},-0.27)
      {$\ell_{\jk}$};
  }

  \foreach \r in {1,...,4}{
    \node[font=\scriptsize]
      at (0.95,{-0.54*(\r+.5)})
      {$V_{\r}$};
  }

  \node[font=\scriptsize]
    at (0.95,{-0.54*5.5})
    {$R_1,R_2,R_3$};

  \foreach \r in {1,...,5}{
    \foreach \c in {1,...,4}{
      \def\g{ab}%

      \ifnum\r=1
        \ifnum\c=1 \def\g{aa}\fi
        \ifnum\c=3 \def\g{bb}\fi
      \fi

      \ifnum\r=2
        \ifnum\c=2 \def\g{aa}\fi
        \ifnum\c=4 \def\g{bb}\fi
      \fi

      \ifnum\r=4
        \ifnum\c=1 \def\g{bb}\fi
        \ifnum\c=2 \def\g{bb}\fi
        \ifnum\c=3 \def\g{aa}\fi
        \ifnum\c=4 \def\g{aa}\fi
      \fi

      \node[font=\scriptsize]
        at ({1.9+(\c-.5)*1.04},
            {-0.54*(\r+.5)})
        {$\g$};
    }
  }

\end{scope}
\end{tikzpicture}%
}
\caption{A four-vertex Hamiltonian-path instance and its
corresponding closed DVI instance.}
\label{fig:hp-dvi-reduction}
\end{figure}

\newcommand{\hpAssignedChain}[6]{%
  \begin{scope}[shift={(#2,0)}]
    \node[hp missing] (#1u1) at (0.70,5.25) {$U_1$};
    \node[hp ref]     (#1r1) at (2.55,5.25) {$R_1$};
    \node[hp missing] (#1u2) at (1.625,3.72) {$U_2$};
    \node[hp ref]     (#1r2) at (3.475,3.72) {$R_2$};
    \node[hp missing] (#1u3) at (2.55,2.19) {$U_3$};
    \node[hp ref]     (#1r3) at (4.40,2.19) {$R_3$};
    \node[hp missing] (#1u4) at (3.475,0.66) {$U_4$};
    \draw (#1u1.east)--(#1r1.west);
    \draw (1.625,5.25)--(#1u2.north);
    \draw (#1u2.east)--(#1r2.west);
    \draw (2.55,3.72)--(#1u3.north);
    \draw (#1u3.east)--(#1r3.west);
    \draw (3.475,2.19)--(#1u4.north);
    \node[font=\scriptsize,anchor=north]
      at ([yshift=-1mm]#1u1.south) {$V_{#3}$};
    \node[font=\scriptsize,anchor=north]
      at ([yshift=-1mm]#1u2.south) {$V_{#4}$};
    \node[font=\scriptsize,anchor=north]
      at ([yshift=-1mm]#1u3.south) {$V_{#5}$};
    \node[font=\scriptsize,anchor=north]
      at ([yshift=-1mm]#1u4.south) {$V_{#6}$};
  \end{scope}%
}

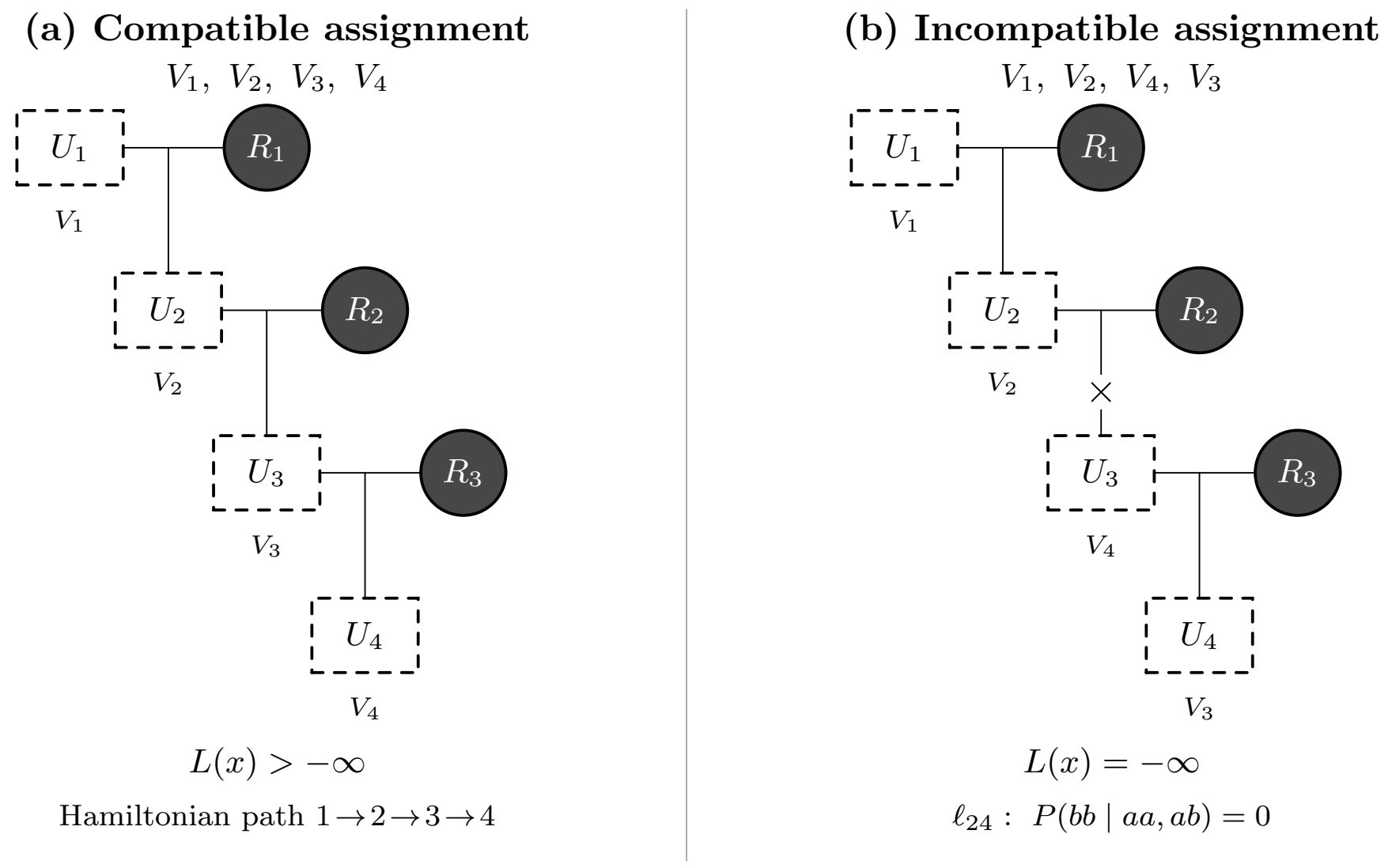
\begin{figure}[H]
\centering
\resizebox{0.7\linewidth}{!}{%
\begin{tikzpicture}[font=\small,line cap=round]
  \node[font=\bfseries] at (2.65,6.35) {(a) Compatible assignment};
  \node at (2.65,5.90) {$V_1,\ V_2,\ V_3,\ V_4$};
  \hpAssignedChain{a}{0}{1}{2}{3}{4}
  \node at (2.65,-0.55) {$L(x)>- \infty$};
  \node[font=\footnotesize] at (2.65,-1.06)
    {Hamiltonian path $1\!\to\!2\!\to\!3\!\to\!4$};

  \draw[black!35] (6.50,6.55)--(6.50,-1.46);

  \node[font=\bfseries] at (10.50,6.35) {(b) Incompatible assignment};
  \node at (10.50,5.90) {$V_1,\ V_2,\ V_4,\ V_3$};
  \hpAssignedChain{b}{7.85}{1}{2}{4}{3}
  \node[fill=white,inner sep=1pt,font=\large] at (10.40,2.95) {$\times$};
  \node at (10.50,-0.55) {$L(x)= - \infty$};
  \node[font=\footnotesize] at (10.50,-1.06)
    {$\ell_{24}:\ P(bb\mid aa,ab)=0$};
\end{tikzpicture}%
}
\caption{A valid assignment corresponds to the Hamiltonian path
  $(v_1,v_2,v_3,v_4)$. An invalid assignment places $V_2$ immediately before $V_4$, violating Mendelian inheritance at locus $\ell_{24}$.}
\label{fig:hp-dvi-assignments}
\end{figure}

\newpage
\section{Linear algebraic representation}\label{sec:linalg}
In this section, we demonstrate how to calculate all likelihoods necessary
for Elston--Stewart, following the reformulation announced in
\Cref{sec:symmetries}. We assume linkage equilibrium, so that it suffices
to work on one locus. Throughout we use the notation of \Cref{methods}:
families $F_j$, slots $I_j$, configurations $s_{i,j}\in\mathcal{V}^{d_j}$,
and weights $w_{i,j}=\log P(s_{i,j}\mid F_j)$.

\subsection{Preliminaries}
We assume $q$ alleles $a_1,\dots,a_q$ with subsequent allele frequencies
$p_1,\dots,p_q$. We write $K=\tfrac12q(q+1)$ for the number of possible
genotypes, which we denote by $\mathcal{A}=\{A_1,\dots,A_K\}$. The
ordering does not matter, but we fix $A_1=a_1a_1$, $A_2=a_1a_2$,
$A_q=a_1a_q$, $A_{q+1}=a_2a_2$ and so on.

Let $G\subseteq\mathcal{A}$ be a collection of genotypes. For a genotype
$g\in\mathcal{A}$, we define the indicator function $\ind_G(g)\in\{0,1\}$
to be $1$ if $g\in G$. Write $2^{\mathcal{A}}$ for the power set of
$\mathcal{A}$. We construct the following maps
$I_F,I_M:2^{\mathcal{A}}\mapsto\{0,1\}^{K^2}$, accompanied by the vectors
$p_F,p_M\in\mathbb{R}^{K^2}$:

\begin{equation}
\bm{I}_F(G)=
\begin{pmatrix}
\ind_{G}(A_1)\\ \ind_{G}(A_1)\\ \vdots\\
\ind_{G}(A_1)\\ \ind_{G}(A_2)\\ \vdots\\
\ind_{G}(A_2)\\ \ind_{G}(A_3)\\ \vdots\\
\vdots\\ \ind_{G}(A_K)
\end{pmatrix},
\quad
\bm{p}_F=
\begin{pmatrix}
p_1^2\\ p_1^2\\ \vdots\\
p_1^2\\ 2p_1p_2\\ \vdots\\
2p_1p_2\\ 2p_1p_3\\ \vdots\\
\vdots\\ p_q^2
\end{pmatrix},
\quad
\bm{I}_M(G)=
\begin{pmatrix}
\ind_{G}(A_1)\\ \ind_{G}(A_2)\\ \vdots\\
\ind_{G}(A_K)\\ \ind_{G}(A_1)\\ \vdots\\
\ind_{G}(A_K)\\ \vdots\\ \ind_{G}(A_1)\\
\vdots\\ \ind_{G}(A_K)
\end{pmatrix},
\quad
\bm{p}_M=
\begin{pmatrix}
p_1^2\\ 2p_1p_2\\ \vdots\\
p_q^2\\ p_1^2\\ \vdots\\
p_q^2\\ \vdots\\ p_1^2\\
\vdots\\ p_q^2
\end{pmatrix}.
\end{equation}
We further write $\bm{F}(G):=\bm{I}_F(G)\odot\bm{p}_F$ and
$\bm{M}(H):=\bm{I}_M(H)\odot\bm{p}_M$, where $\odot$ denotes the
element-wise Hadamard product.
Assuming no inbreeding, the genotypes $g_F$ and $g_M$ of a father $F$ and
a mother $M$ are independent. Now observe that the joint likelihood that
$g_F\in G$ and $g_M\in H$ can be expressed as
\begin{equation}
\langle\bm{F}(G),\bm{M}(H)\rangle.
\end{equation}
In practice, $G=\{g\}$ or $G=\mathcal{A}$ since individuals are usually
typed or unknown. The same holds for $H$.

The genotypes of the children $C_1,\dots,C_N$ of $F$ and $M$ are of course
dependent on $F$ and $M$. For them we need a different notation. We start
by writing
\begin{equation}
\bm{I}_C(G)=
\begin{pmatrix}
\ind_{G}(A_1)\\ \vdots\\ \ind_{G}(A_K)
\end{pmatrix}.
\end{equation}

The genotypes of children are independent conditionally on their parents'
genotypes. Therefore, let $C$ be any child of $F$ and $M$. For each
$A_i\in\mathcal{A}$, we define
\begin{equation}
\bm{A}_i=
\begin{pmatrix}
P(C=A_i\mid F=A_1,M=A_1)\\
P(C=A_i\mid F=A_1,M=A_2)\\
\vdots\\
P(C=A_i\mid F=A_2,M=A_K)\\
P(C=A_i\mid F=A_2,M=A_1)\\
\vdots\\
P(C=A_i\mid F=A_K,M=A_K)
\end{pmatrix}.
\end{equation}
These vectors are stored in a $K^2\times K$ matrix as follows:
\begin{equation}
\bm{P}_C:=
\begin{pmatrix}
\mid & \mid & & \mid\\
\bm{A}_1 & \bm{A}_2 & \hldots & \bm{A}_K\\
\mid & \mid & & \mid\\
\end{pmatrix}.
\end{equation}
As with $\bm{F}$ and $\bm{M}$, we write $\bm{C}(G):=\bm{P}_C\bm{I}_C(G)$.

Now, let $\mathcal{P}$ denote the pedigree where $C$ is a child of $F$ and
$M$. Again assuming no inbreeding, and writing $g_C$ for the genotype of
$C$, we have
\begin{equation}
P(g_F,g_M,g_C\mid\mathcal{P})
=\langle\bm{C}(g_C),\bm{F}(g_F)\odot\bm{M}(g_M)\rangle,
\end{equation}
with the slightly abusive notation $g_F=\{g_F\}$.

If $\mathcal{P}$ contains multiple children $C_1,\dots,C_N$ with respective
genotypes $g_{C_1},\dots,g_{C_N}$, the resulting expression is
surprisingly simple due to conditional independence:
\begin{equation}
P(g_F,g_M,g_{C_1},\dots,g_{C_N}\mid\mathcal{P})
=\Bigl\langle\bigodot_{j=1}^N\bm{C}(g_{C_j}),\,
\bm{F}(g_F)\odot\bm{M}(g_M)\Bigr\rangle.
\end{equation}

The Hadamard product is associative with respect to the inner product, so
that we also have
\begin{equation}
\begin{split}
P(g_F,g_M,g_{C_1},\dots,g_{C_N}\mid\mathcal{P})
&=\Bigl\langle\bigodot_{j=1}^N\bm{C}(g_{C_j})\odot\bm{M}(g_M),\,
\bm{F}(g_F)\Bigr\rangle\\
&=\Bigl\langle\bigodot_{j=1}^N\bm{C}(g_{C_j})\odot\bm{F}(g_F),\,
\bm{M}(g_M)\Bigr\rangle,
\end{split}
\end{equation}
which will be useful later on for transmitting messages. In the nuclear
family example of \Cref{sec:symmetries}, $\bm{v}=\bm{F}\odot\bm{M}$ and
$\bm{w}=\bm{C}$ with symbolic indicators, and the inner product
\eqref{nuclinprodr} is exactly this identity.

\subsection{Sending a message}
When messages are sent, we again make use of conditional independence.
Since there are no loops, dependencies between peels always go through
exactly one individual. If we fix the profile of this individual, the
separate parts this individual connects are independent.

\subsubsection{Parent to parent}
If a parent sends a message, this is one side of the inner product without
the parent that is sending. Let $\mathcal{P}$ be the pedigree a parent is
sending from. Then we write
\begin{equation}
\begin{split}
\bm{T}_F&:=\bigodot_{j=1}^N\bm{C}(g_{C_j})\odot\bm{M}(g_M),\\
\bm{T}_M&:=\bigodot_{j=1}^N\bm{C}(g_{C_j})\odot\bm{F}(g_F).
\end{split}
\end{equation}

Fix some profile $A_i$. We write $P(\mathcal{P}\backslash X)$ for the
likelihood of the given profiles in the pedigree, but without individual
$X$. Now observe that
\begin{equation}
\begin{split}
P(\mathcal{P}\backslash F\mid F=A_i)
&=\langle\bm{T}_F,\bm{I}_F(A_i)\rangle,\\
P(\mathcal{P}\backslash M\mid M=A_i)
&=\langle\bm{T}_M,\bm{I}_M(A_i)\rangle.
\end{split}
\end{equation}

This motivates the following definition for the paternal resp.\ maternal
\textit{sideways message}:
\begin{equation}
\bm{S}_F=
\begin{pmatrix}
\langle\bm{T}_F,\bm{I}_F(A_1)\rangle\\
\langle\bm{T}_F,\bm{I}_F(A_1)\rangle\\ \vdots\\
\langle\bm{T}_F,\bm{I}_F(A_1)\rangle\\
\langle\bm{T}_F,\bm{I}_F(A_2)\rangle\\ \vdots\\
\langle\bm{T}_F,\bm{I}_F(A_2)\rangle\\
\langle\bm{T}_F,\bm{I}_F(A_3)\rangle\\ \vdots\\
\vdots\\
\langle\bm{T}_F,\bm{I}_F(A_K)\rangle
\end{pmatrix},
\quad
\bm{S}_M=
\begin{pmatrix}
\langle\bm{T}_M,\bm{I}_M(A_1)\rangle\\
\langle\bm{T}_M,\bm{I}_M(A_2)\rangle\\ \vdots\\
\langle\bm{T}_M,\bm{I}_M(A_K)\rangle\\
\langle\bm{T}_M,\bm{I}_M(A_1)\rangle\\ \vdots\\
\langle\bm{T}_M,\bm{I}_M(A_K)\rangle\\ \vdots\\
\langle\bm{T}_M,\bm{I}_M(A_1)\rangle\\
\vdots\\
\langle\bm{T}_M,\bm{I}_M(A_K)\rangle
\end{pmatrix}.
\end{equation}

Suppose that $F$ is a father in pedigrees $\mathcal{P}_1$ and
$\mathcal{P}_2$. Let $\bm{S}_F$ be calculated on $\mathcal{P}_1$. Then by
conditional independence,
\begin{equation}
P(\mathcal{P}_1,\mathcal{P}_2)
=\sum_F P(\mathcal{P}_1\backslash F\mid F)\cdot
P(\mathcal{P}_2\backslash F\mid F)\cdot P(F).
\end{equation}
This is exactly what the inner product does:
\begin{equation}
\begin{split}
\Bigl\langle\bigodot_{j=1}^N\bm{C}(g_{C_j}),\,
\bm{F}(g_F)\odot\bm{S}_F\odot\bm{M}(g_M)\Bigr\rangle
&=\Bigl\langle\bigodot_{j=1}^N\bm{C}(g_{C_j})\odot\bm{M}(g_M),\,
\bm{F}(g_F)\odot\bm{S}_F\Bigr\rangle\\
&=\Bigl\langle\bigodot_{j=1}^N\bm{C}(g_{C_j})\odot\bm{M}(g_M),\,
\bm{p}_F\odot\bm{S}_F\Bigr\rangle.
\end{split}
\end{equation}
The right side contains the entries
$P(F=A_i)\cdot P(\mathcal{P}_1\backslash F\mid F=A_i)
=P(F=A_i,\mathcal{P}_1\backslash F)$, whereas the inner product with the
left side induces multiplication with
$P(\mathcal{P}_2\backslash F\mid F=A_i)$.

\subsubsection{Parent to child}
An untyped parent can also be a child in the pedigree it sends the
message to. In that case, we need to ``wrap'' the message in a matrix
since that is the way we represent children. The paternal resp.\ maternal
\textit{upward message} is defined as
\begin{equation}
\begin{split}
\bm{U}_F&:=
\begin{pmatrix}
\langle\bm{T}_F,\bm{I}_F(A_1)\rangle &
\langle\bm{T}_F,\bm{I}_F(A_2)\rangle &
\hldots &
\langle\bm{T}_F,\bm{I}_F(A_K)\rangle\\
\vdots & \vdots & & \vdots\\
\langle\bm{T}_F,\bm{I}_F(A_1)\rangle &
\langle\bm{T}_F,\bm{I}_F(A_2)\rangle &
\hldots &
\langle\bm{T}_F,\bm{I}_F(A_K)\rangle
\end{pmatrix},\\[4pt]
\bm{U}_M&:=
\begin{pmatrix}
\langle\bm{T}_M,\bm{I}_M(A_1)\rangle &
\langle\bm{T}_M,\bm{I}_M(A_2)\rangle &
\hldots &
\langle\bm{T}_M,\bm{I}_M(A_K)\rangle\\
\vdots & \vdots & & \vdots\\
\langle\bm{T}_M,\bm{I}_M(A_1)\rangle &
\langle\bm{T}_M,\bm{I}_M(A_2)\rangle &
\hldots &
\langle\bm{T}_M,\bm{I}_M(A_K)\rangle
\end{pmatrix}.
\end{split}
\end{equation}
The idea is exactly the same as for sideways messages. The probability of
the two pedigrees is again conditionally independent:
\begin{equation}
\begin{split}
\langle\bm{T}_M,\bm{I}_M(A_i)\rangle\cdot P(C=A_i\mid F,M)
&=P(\mathcal{P}\backslash C\mid C=A_i)\cdot P(C=A_i\mid F,M)\\
&=P(\mathcal{P}\backslash C\mid C=A_i,F,M)\cdot P(C=A_i\mid F,M)\\
&=P(C,\mathcal{P}\backslash C\mid F,M).
\end{split}
\end{equation}

\subsubsection{Child to parent}
Every peel is sending at most one message. Without loss of generality, we
assume that if a child sends the message, it is $C_1$ that is sending.
We first write
\begin{equation}
\bm{T}_C:=\bigodot_{j=2}^N\bm{C}(g_{C_j})\odot\bm{M}(g_M)\odot\bm{F}(g_F).
\end{equation}
Now we define the paternal resp.\ maternal \textit{downward message} as
follows:
\begin{equation}
\bm{D}_F=
\begin{pmatrix}
\langle\bm{T}_C,\bm{C}(A_1)\rangle\\
\langle\bm{T}_C,\bm{C}(A_1)\rangle\\ \vdots\\
\langle\bm{T}_C,\bm{C}(A_1)\rangle\\
\langle\bm{T}_C,\bm{C}(A_2)\rangle\\ \vdots\\
\langle\bm{T}_C,\bm{C}(A_2)\rangle\\
\langle\bm{T}_C,\bm{C}(A_3)\rangle\\ \vdots\\
\vdots\\
\langle\bm{T}_C,\bm{C}(A_K)\rangle
\end{pmatrix},
\quad
\bm{D}_M=
\begin{pmatrix}
\langle\bm{T}_C,\bm{C}(A_1)\rangle\\
\langle\bm{T}_C,\bm{C}(A_2)\rangle\\ \vdots\\
\langle\bm{T}_C,\bm{C}(A_K)\rangle\\
\langle\bm{T}_C,\bm{C}(A_1)\rangle\\ \vdots\\
\langle\bm{T}_C,\bm{C}(A_K)\rangle\\ \vdots\\
\langle\bm{T}_C,\bm{C}(A_1)\rangle\\
\vdots\\
\langle\bm{T}_C,\bm{C}(A_K)\rangle
\end{pmatrix}.
\end{equation}
Observe that
\begin{equation}
P(\mathcal{P}\mid C_1=A_i)
=\langle\bm{T}_C,\bm{P}_C\bm{I}_C(A_i)\rangle
=\langle\bm{T}_C,\bm{C}(A_i)\rangle.
\end{equation}
This message is different in the sense that it is already multiplied with
$P(C=A_i)$, because this is conditional even within its peel. Thus, when
this message is sent, the parent receiving this message should not also
multiply this with $\bm{p}$. (The good news is that a child can only be a
child once.)

\subsection{Receiving a message}
When an individual receives a message, we assume it is always unknown. If
it is an upward or sideways message, we simply take the Hadamard product
$\bm{S}_F\odot\bm{p}_F$. For an upward message into a child role, it is
$\bm{U}_F$ composed with $\bm{C}$. For a downward message, $\bm{F}$ is
replaced by $\bm{D}_F$, as it does not need the Hadamard product with
$\bm{p}_F$.

A peel can receive multiple messages, but sends only one. In the ``final''
peel, messages can only be received. Assuming all members receive, we can
have for instance
\begin{equation}
\Bigl\langle\bigodot_{j=1}^N\bm{C}(g_{C_j})\odot\bm{U}_j,\,
\bm{S}_F\odot\bm{p}_F\odot\bm{D}_M\Bigr\rangle
\end{equation}
for the final peel. Here, the mother received a message as a child, the
father had other children, and the $\bm{U}_j$ denote the messages of the
children they transmit as parents. If this is not the final peel, a
subsequent $\bm{T}$ can be constructed by selecting the appropriate
factors.

Parents can receive multiple messages $\bm{T}^i_F$ if they have children
with several partners. In such cases, one takes the Hadamard product of
these messages, and sends off a complete message
$\bm{T}_F=\bigodot_i\bm{T}^i_F$.

As mentioned before, kids only have one pair of parents. However, if they
need to send a downward message, while also having children with multiple
partners, the complete message is
$\bm{T}_C\odot\bigodot_i\bm{T}^i_F$. Since there is only one Hadamard
factor $\bm{T}_C$, it is no problem that this already contains $\bm{p}$.

\subsubsection{Example of a sideways message}
A typical example is a half-sibling pedigree:

\begin{figure}[H]
\centering
\begin{tikzpicture}[
  scale=0.8,
  every node/.style={transform shape},
  ped male/.style   = {draw, rectangle, minimum size=8.5mm, inner sep=0pt},
  ped female/.style = {draw, circle,    minimum size=8.5mm, inner sep=0pt},
  ped ref/.style    = {fill=black!72, text=white, thick},
  ped miss/.style   = {dashed, thick},
  ped unk/.style    = {fill=black!7}
]
\node[ped female, ped miss] (M1) at (-3,0) {\scriptsize $M_1$};
\node[ped male,   ped unk]  (F)  at ( 0,0) {\scriptsize $F$};
\node[ped female, ped miss] (M2) at ( 3,0) {\scriptsize $M_2$};
\node[ped male, ped ref] (C1) at (-1.5,-2.8) {\scriptsize $C_1$};
\node[ped male, ped ref] (C2) at ( 1.5,-2.8) {\scriptsize $C_2$};
\draw (M1.east) -- (F.west);
\draw (F.east)  -- (M2.west);
\draw (-1.5,0) -- (C1.north);
\draw ( 1.5,0) -- (C2.north);
\node[anchor=north] at ([yshift=-1.5mm]C1.south) {\scriptsize $A/A$};
\node[anchor=north] at ([yshift=-1.5mm]C2.south) {\scriptsize $A/a$};
\end{tikzpicture}
\caption{Pedigree with two missing mothers, an untyped father,
and two typed children with genotypes $A/A$ and $A/a$.}
\label{fig:two-mothers-pedigree}
\end{figure}
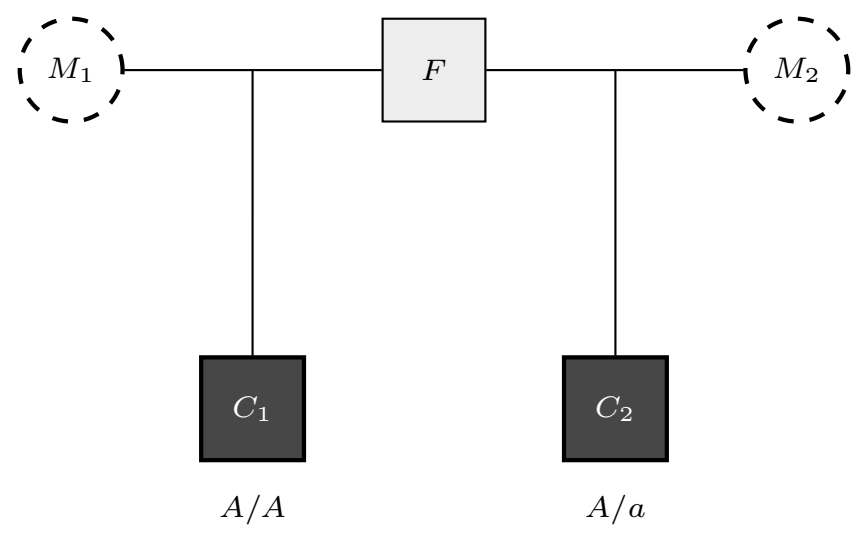

If we assume the message via $F$ goes from left to right, we have
\begin{equation}
\begin{split}
\bm{T}_F
&=\bm{C}(AA)\odot\bm{p}_M
=
\begin{pmatrix}
0\\0\\0\\0\\\tfrac14\\\tfrac12\\0\\\tfrac12\\1
\end{pmatrix}
\odot
\begin{pmatrix}
\tfrac14\\\tfrac12\\\tfrac14\\
\tfrac14\\\tfrac12\\\tfrac14\\
\tfrac14\\\tfrac12\\\tfrac14
\end{pmatrix}
=
\begin{pmatrix}
0\\0\\0\\0\\\tfrac18\\\tfrac18\\0\\\tfrac14\\\tfrac14
\end{pmatrix}.
\end{split}
\end{equation}
Now we calculate
\begin{equation}
\langle\bm{T}_F,\bm{I}_F(A_1)\rangle=0,\quad
\langle\bm{T}_F,\bm{I}_F(A_2)\rangle=\tfrac14,\quad
\langle\bm{T}_F,\bm{I}_F(A_3)\rangle=\tfrac12,
\end{equation}
so that
\begin{equation}
\bm{S}_F=
\begin{pmatrix}
0\\0\\0\\\tfrac14\\\tfrac14\\\tfrac14\\\tfrac12\\\tfrac12\\\tfrac12
\end{pmatrix}.
\end{equation}
This is the message the father receives through the left. Now we calculate
\begin{equation}
\begin{split}
LL
&=\Bigl\langle\bigodot_{j=1}^N\bm{C}(g_{C_j})\odot\bm{U},\,
\bm{S}_F\odot\bm{p}_F\odot\bm{D}_M\Bigr\rangle\\
&=\langle\bm{C}(Aa),\,\bm{S}_F\odot\bm{p}_F\odot\bm{M}\rangle\\
&=\langle\bm{C}(Aa),\,\bm{S}_F\odot\bm{p}_F\odot\bm{p}_M\rangle\\
&=\Biggl\langle
\begin{pmatrix}
0\\\tfrac12\\1\\\tfrac12\\\tfrac12\\\tfrac12\\1\\\tfrac12\\0
\end{pmatrix},
\begin{pmatrix}
0\\0\\0\\\tfrac1{32}\\\tfrac1{16}\\\tfrac1{32}\\\tfrac1{32}\\\tfrac1{16}\\\tfrac1{32}
\end{pmatrix}
\Biggr\rangle
=\tfrac18.
\end{split}
\end{equation}

\end{document}